\documentclass[reqno,twoside,11pt]{amsart}
 
\hoffset - 1.5cm

\usepackage{graphicx}
\usepackage{amstext}
\usepackage[T1]{fontenc}
\usepackage{pstricks,pst-node,pst-text,pst-3d}
\usepackage{color}
\usepackage{amsmath}
\usepackage{amsthm}
\usepackage{amssymb}
\usepackage{eufrak}
\usepackage[mathscr]{euscript}

\newtheorem{Theorem}{Theorem}[section]
\newtheorem{Fact}{Fact}[section]
\newtheorem{Lemma}{Lemma}[section]
\newtheorem{Proposition}{Proposition}[section]
\newtheorem{Corollary}{Corollary}[section]
\theoremstyle{definition}
\newtheorem{Definition}{Definition}[section]
\newtheorem{Example}{Example}[section]
\newtheorem{Remark}{Remark}[section]

\newcommand{\numsec}{\setcounter{Theorem}{0}\setcounter{Definition}{0}
\setcounter{Remark}{0} \setcounter{Lemma}{0} \setcounter{Fact}{0}
\setcounter{Proposition}{0} \setcounter{Corollary}{0}
\setcounter{Example}{0} \setcounter{equation}{0}
\setcounter{Property}{0}\renewcommand\theequation{\arabic{section}.\arabic{equation}}
\renewcommand\theTheorem{\arabic{section}.\arabic{Theorem}}
\renewcommand\theDefinition{\arabic{section}.\arabic{Definition}}
\renewcommand\theRemark{\arabic{section}.\arabic{Remark}}
\renewcommand\theLemma{\arabic{section}.\arabic{Lemma}}
\renewcommand\theFact{\arabic{section}.\arabic{Fact}}
\renewcommand\theProposition{\arabic{section}.\arabic{Proposition}}
\renewcommand\theCorollary{\arabic{section}.\arabic{Corollary}}
\renewcommand\theExample{\arabic{section}.\arabic{Example}}
\renewcommand\theProperty{\arabic{section}.\arabic{Property}}}
\newcommand{\numsubsec}{\setcounter{Theorem}{0}\setcounter{Definition}{0}
\setcounter{Remark}{0} \setcounter{Lemma}{0} \setcounter{Fact}{0}
\setcounter{Proposition}{0} \setcounter{Corollary}{0}
\setcounter{Example}{0} \setcounter{equation}{0}
\setcounter{Property}{0}\renewcommand\theequation{\arabic{section}.\arabic{subsection}.\arabic{equation}}
\renewcommand\theTheorem{\arabic{section}.\arabic{subsection}.\arabic{Theorem}}
\renewcommand\theDefinition{\arabic{section}.\arabic{subsection}.\arabic{Definition}}
\renewcommand\theRemark{\arabic{section}.\arabic{subsection}.\arabic{Remark}}
\renewcommand\theLemma{\arabic{section}.\arabic{subsection}.\arabic{Lemma}}
\renewcommand\theFact{\arabic{section}.\arabic{subsection}.\arabic{Fact}}
\renewcommand\theProposition{\arabic{section}.\arabic{subsection}.\arabic{Proposition}}
\renewcommand\theCorollary{\arabic{section}.\arabic{subsection}.\arabic{Corollary}}
\renewcommand\theExample{\arabic{section}.\arabic{subsection}.\arabic{Example}}
\renewcommand\theProperty{\arabic{section}.\arabic{subsection}.\arabic{Property}}}

\newcommand{\numsubsubsec}{\setcounter{Theorem}{0}\setcounter{Definition}{0}
\setcounter{Remark}{0} \setcounter{Lemma}{0} \setcounter{Fact}{0}
\setcounter{Proposition}{0} \setcounter{Corollary}{0}
\setcounter{Example}{0} \setcounter{equation}{0}
\setcounter{Property}{0}\renewcommand\theequation{\arabic{section}.\arabic{subsection}.\arabic{subsubsection}.\arabic{equation}}
\renewcommand\theTheorem{\arabic{section}.\arabic{subsection}.\arabic{subsubsection}.\arabic{Theorem}}
\renewcommand\theDefinition{\arabic{section}.\arabic{subsection}.\arabic{subsubsection}.\arabic{Definition}}
\renewcommand\theRemark{\arabic{section}.\arabic{subsection}.\arabic{subsubsection}.\arabic{Remark}}
\renewcommand\theLemma{\arabic{section}.\arabic{subsection}.\arabic{subsubsection}.\arabic{Lemma}}
\renewcommand\theFact{\arabic{section}.\arabic{subsection}.\arabic{subsubsection}.\arabic{Fact}}
\renewcommand\theProposition{\arabic{section}.\arabic{subsection}.\arabic{subsubsection}.\arabic{Proposition}}
\renewcommand\theCorollary{\arabic{section}.\arabic{subsection}.\arabic{subsubsection}.\arabic{Corollary}}
\renewcommand\theExample{\arabic{section}.\arabic{subsection}.\arabic{subsubsection}.\arabic{Example}}
\renewcommand\theProperty{\arabic{section}.\arabic{subsection}.\arabic{subsubsection}.\arabic{Property}}}

\newcommand{\ba}{\begin{array}}
\newcommand{\bc}{\begin{center}}
\newcommand{\bd}{\begin{description}}
\newcommand{\bdm}{\begin{displaymath}}
\newcommand{\be}{\begin{enumerate}}
\newcommand{\beq}{\begin{equation}}
\newcommand{\bdf}{\begin{Definition}}
\newcommand{\bex}{\begin{Example}}
\newcommand{\bft}{\begin{Fact}}
\newcommand{\bl}{\begin{Lemma}}
\newcommand{\bp}{\begin{Proposition}}
\newcommand{\br}{\begin{Remark}}
\newcommand{\bt}{\begin{Theorem}}
\newcommand{\bco}{\begin{Corollary}}
\newcommand{\bh}{\begin{Hipothesis}}
\newcommand{\ea}{\end{array}}
\newcommand{\ec}{\end{center}}
\newcommand{\ed}{\end{description}}
\newcommand{\edm}{\end{displaymath}}
\newcommand{\ee}{\end{enumerate}}
\newcommand{\eeq}{\end{equation}}
\newcommand{\edf}{\end{Definition}}
\newcommand{\eex}{\end{Example}}
\newcommand{\eft}{\end{Fact}}
\newcommand{\el}{\end{Lemma}}
\newcommand{\ep}{\end{Proposition}}
\newcommand{\er}{\end{Remark}}
\newcommand{\et}{\end{Theorem}}
\newcommand{\eco}{\end{Corollary}}
\newcommand{\eh}{\end{Hipothesis}}

\newcommand{\bA}{\mathbb{A}}

\newcommand{\bC}{\mathbb{C}}

\newcommand{\bH}{\mathbb{H}}
\newcommand{\bI}{\mathbb{I}}

\newcommand{\bN}{\mathbb{N}}

\newcommand{\bR}{\mathbb{R}}

\newcommand{\bV}{\mathbb{V}}
\newcommand{\bW}{\mathbb{W}}
\newcommand{\bX}{\mathbb{X}}
\newcommand{\bY}{\mathbb{Y}}
\newcommand{\bZ}{\mathbb{Z}}

\newcommand{\cC}{\mathcal{C}}

\newcommand{\cE}{\mathcal{E}}
\newcommand{\cF}{\mathcal{F}}

\newcommand{\cH}{\mathcal{H}}
\newcommand{\cI}{\mathcal{I}}

\newcommand{\cL}{\mathcal{L}}

\newcommand{\cN}{\mathcal{N}}
\newcommand{\cO}{\mathcal{O}}

\newcommand{\cR}{\mathcal{R}}

\newcommand{\cT}{\mathcal{T}}
\newcommand{\cU}{\mathcal{U}}

\newcommand{\im}{\mathrm{ im \;}}

\numberwithin{equation}{section} \errorcontextlines=0

\newcommand{\vp}{\varphi}

\newcommand{\lin}{\mathrm{span}}

\newcommand{\morse}{\mathrm{m^-}}

\newcommand{\sone}{S^1}
\newcommand{\sotwo}{SO(2)}

\newcommand{\ds}{\displaystyle}
\newcommand{\nt}{\noindent}

\newcommand{\on}{O(n,\bR)}

\newcommand{\subg}{\overline{\mathrm{sub}}(G)}
\newcommand{\subgc}{\overline{\mathrm{sub}}[G]}
\newcommand{\subh}{\overline{\mathrm{sub}}(H)}

\newcommand{\gcw}{G\text{-CW}}
\newcommand{\gls}{$G$-$\mathcal{LS}$}

\newcommand{\chig}{\chi_G}
\newcommand{\chih}{\chi_H}

\newcommand{\upsg}{\Upsilon_G}
\newcommand{\upsh}{\Upsilon_H}
\newcommand{\upssone}{\Upsilon_{\sone}}
\newcommand{\ci}{\cC\cI}
\newcommand{\cig}{\ci_G}
\newcommand{\cih}{\ci_H}

\newcommand{\scig}{\mathscr{CI}_G}
\newcommand{\sci}{\mathscr{CI}}
\newcommand{\lpm}{\lambda_{\pm}}

\newcommand{\tpsil}{\widetilde{\Psi}_{\lambda}}
\newcommand{\tangent}{T_{x_0}^{\perp} G(x_0)}

\newcommand{\h}{\mathbb{H}}
\newcommand{\m}{\text{mult}}

\begin{document}

\title[Symmetric Liapunov center theorem for minimal orbit]{Symmetric Liapunov center theorem for minimal orbit}

\author{Ernesto P\'{e}rez-Chavela$^{1)}$}
\address{$^{1)}$ Departamento de Matem\'{a}ticas \\ Instituto Tecnologico Aut\'{o}nomo de M\'{e}xico (ITAM), R\'{o} Hondo 1, Col. Progreso Tizap\'{a}n
M\'{e}xico D.F. 01080 \\ M\'{e}xico}

\author{S{\l}awomir Rybicki$^{2)}$}
\address{$^{2)}$Faculty of Mathematics and Computer Science\\
Nicolaus Copernicus University \\
PL-87-100 Toru\'{n} \\ ul. Chopina $12 \slash 18$ \\
Poland}

\author{Daniel Strzelecki$^{2)}$}

\email{ernesto.perez@itam.mx (E. P\'{e}rez-Chavela)}
\email{rybicki@mat.umk.pl (S. Rybicki)}
\email{danio@mat.umk.pl (D. Strzelecki)}

\date{\today}


\begin{abstract} Using the techniques of equivariant bifurcation theory  
we prove the existence of non-stationary periodic solutions of $\Gamma$-symmetric  systems  $\ddot q(t)=-\nabla U(q(t))$ in any neighborhood of an isolated orbit of minima $\Gamma(q_0)$ of the potential $U.$ We show the strength of our result by proving the existence of new families of periodic orbits in the Lennard-Jones two- and three-body problems and in the Schwarzschild three-body problem. 

\end{abstract}

\maketitle

%%%%%%%%%%%%%%%%%%%%%%%%%%%%%%%%%%%%%%%%%%%%%%%%%%%%%%%%%%%%%%
\section{Introduction}
\numsec

The study of the existence of non-stationary periodic solutions of autonomous ordinary differential equations has a long history. Particular attention was paid to the study of the existence of such solutions in a neighborhood  of isolated equilibria, see for instance \cite{FUGORY, HENRARD, MARPAR,  MAWI, MOSER, SCHMIDT} and references therein. Of course this list is far from being complete.

One of the most famous theorems concerning the existence of periodic solutions of ordinary differential equations is the celebrated Liapunov center theorem. 

Consider a second order  system $\ddot q(t)  =  - \nabla U(q(t)),$ where $U \in C^2(\bR^n,\bR), \nabla U(0)=0$ and $ \det \nabla^2 U(0) \neq 0.$ Let $\sigma(\nabla^2 U(0))$ be the spectrum of the Hessian  $\nabla^2 U(0).$ The Liapunov center theorem says that if  $\sigma(\nabla^2 U(0)) \cap (0,+\infty) = \{\beta_1^2,\ldots, \beta_m^2\}$ for $\beta_1 >\ldots > \beta_m > 0$ then for  $\beta_{j_0}$ satisfying $\beta_1 \slash \beta_{j_0}, \ldots, \beta_{j _0-1} \slash \beta_{j_0} \not \in \bN,$   there is a sequence $\{q_k(t)\}$ of periodic solutions of the system  
\beq \label{newsys} \ddot q(t)  =  - \nabla U(q(t)),\eeq 
with amplitude  tending to zero and the  minimal period tending to $2 \pi \slash \beta_{j_0}.$ The proof of this theorem  can be found in \cite{MAWI}, see also \cite{BERGER1,BERGER2}. 
Generalizations of the Liapunov center theorem  were developed in many directions. In \cite{BERGER1, BERGER2, DARY,  PCHRYST, SZULKIN} one can find some of them.

Let  $\Omega \subset \bR^n$ be an open and $\Gamma$-invariant subset of $\bR^n$ considered as a representation of a compact Lie group $\Gamma$. Assume that $q_0 \in \Omega$ is a critical point of the $\Gamma$-invariant potential $U : \Omega \to \bR$  of class $C^2$. Since for all $\gamma \in \Gamma$ the equality $U(\gamma  q_0)=U(q_0)$ holds and $\nabla U(q_0)=0$, the orbit 
$\Gamma(q_0)=\{\gamma q_0 : \gamma \in \Gamma\}$ consists of critical points of $U$ i.e. $\Gamma(q_0) \subset (\nabla U)^{-1}(0).$ Note that if $\dim \Gamma \geqq 1$ then it can happen that $\dim \Gamma(q_0) \geqq 1$ i.e. the critical point $q_0$ is not isolated in $(\nabla U)^{-1}(0)$. That is why for higher-dimensional orbits $\Gamma(q_0)$ we can not apply the classical Liapunov center theorem.

 In \cite{PCHRYST} we have proved the Symmetric Liapunov center theorem for non-degenerate orbit of critical points $\Gamma(q_0)$ i.e.  we have assumed that $\dim \Gamma(q_0) = \dim \ker \nabla^2 U(q_0).$ More precisely, with the additional hypothesis that the isotropy group  $\Gamma_{q_0}=\{\gamma \in \Gamma : \gamma q_0=q_0\}$ is trivial and that there is at least one positive eigenvalue of the Hessian $\nabla^2 U(q_0),$ we have proved the existence of non-stationary periodic solutions of system \eqref{newsys} in any neighborhood of the orbit $\Gamma(q_0).$ Moreover, we are able to control the minimal period of these solutions in terms of the positive eigenvalues of $\nabla^2 U(q_0).$

For the Lennard-Jones and Schwarzschild problems discussed in the last section  there are isolated degenerate circles ($\Gamma=\sotwo$-orbits)  of stationary solutions  which consist of minima of  the corresponding potentials. We underline that we are not able to study the non-stationary periodic solutions of these problems applying the classical Liapunov center theorem because these equilibria are not isolated.
We also emphasize that since these orbits are degenerate, we either can not study the non-stationary periodic solutions of these problems applying the Symmetric Liapunov center theorem for non-degenerate orbit proved in  \cite{PCHRYST}.
Therefore there is a natural need to prove the Symmetric Liapunov center theorem for isolated orbits of minima.

The inspiration for writing this article, in addition to the discussion above,  was a nice paper  of Rabinowitz \cite{RAB4}, where the author proved that the Brouwer index of an isolated minimum of a potential of the class $C^1$ is equal to $1.$ This result  was also proved later by Amann  \cite{AMANN}. 
 
The goal of this paper is to prove the Symmetric Liapunov center theorem for  an isolated orbit of minima of the potential $U$.
Our main result  is the following.

\bt\label{main-theo} [Symmetric Liapunov center theorem for a minimal orbit]  Let $U:\Omega\to\bR$ be a $\Gamma$-invariant potential of the class $C^2$ and $q_0 \in \Omega$. Assume that
\begin{enumerate}
\item the orbit $\Gamma(q_0)$ consists of minima of potential $U$,
\item the orbit  $\Gamma(q_0)$ is isolated in $(\nabla U)^{-1}(0)$,
\item the isotropy group $\Gamma_{q_0}$ is trivial,
\item $\sigma(\nabla^2 U(q_0)) \cap (0,+\infty) = \{\beta_1^2,\ldots, \beta_m^2\}$, $\beta_1>\beta_2>\ldots>\beta_m>0$ and $m\geq 1$.
\end{enumerate}
Then for any $\beta_{j_0}$ such that $\beta_j \slash \beta_{j_0}\not \in \bN$ for all $j\neq j_0$ there exists a sequence $(q_k(t))$ of periodic solutions of the system $\ddot q(t)=-\nabla U(q(t))$ with a sequence of minimal periods $(T_k)$ such that $\mathrm{dist} (\Gamma(q_0),q_k([0,T_k]))\to 0$ and $T_k\to 2 \pi \slash \beta_{j_0}$ as $k\to\infty$.
\et
To prove the above  theorem we apply techniques of the $(\Gamma \times \sone)$-equivariant bifurcation theory. 
We present the problem of  the existence of periodic solutions of  system  \eqref{newsys} as $(\Gamma \times \sone)$-symmetric variational bifurcation problem i.e. we  look for periodic solutions of system \eqref{newsys} as a $(\Gamma \times \sone)$-orbits of critical points of a family  $(\Gamma \times \sone)$-invariant   functionals defined on a suitably chosen orthogonal Hilbert representation of $\Gamma \times \sone$.
As topological tools we apply the $(\Gamma \times \sone)$-equivariant Conley index due to Izydorek, see \cite{IZYDOREK}, and the degree for $(\Gamma \times \sone)$-equivariant gradient operators due to  Go{\l}\c{e}biewska and the second author, see \cite{GORY}.

More precisely, we have proved  changes of the equivariant Conley index and the degree for equivariant gradient operators along the  family $\Gamma(q_0) \times (0,+\infty) \subset \h^1_{2\pi}  \times (0,+\infty) $ of stationary solutions of the following  system 
\beq \label{newsysp}
\left\{\ba{rcl} \ddot q(t) & = & - \lambda^2 \nabla U(q(t)), \\ q(0) & = & q(2\pi), \\ \dot q(0) & = & \dot q(2\pi).  \ea  \right.
\eeq
 We emphasize that change of the Conley index implies the existence of a local bifurcation of periodic solutions of system \eqref{newsysp}, whereas a change of the degree implies the existence of a global bifurcation of periodic solutions  of system \eqref{newsysp} satisfying the Rabinowitz type alternative.

In order to get an accurate model to study the action of the intermolecular
and gravitational forces at the same time, many authors from physics,
astrophysics, astronomy, cosmology  and chemistry have introduced new kinds of
potentials, with a structure different from the classical Newtonians and
Coulombians potentials. In this way, potentials that have been used very
often in those branches of the science are the Lennard-Jones  and the Schwarzschild potentials.  

In the last section we apply Theorem \ref{main-theo} to the study non-stationary periodic solutions of the Lennard-Jones and the Schwarzschild $N$-body problems.

After introduction our paper is organized as follows. 
In Section \ref{preliminaria}  we summarize without proofs the relevant material on equivariant topology  and prove some preliminary results. Throughout this section $G$ stands for a compact Lie group.  Since admissible pairs play a crucial role in our reasonings, the notion of an admissible pair is given in Definition \ref{admissible}.  In Definition \ref{G-spectrum} we introduce the notion of $G$-equivariant spectrum. $G$-equivariant Euler characteristic of $G$-spectrum is given by formula \eqref{gchargspe}.  The $G$-equivariant Conley index which we apply in this article is a $G$-homotopy type of a $G$-spectrum. In Theorem \ref{cio} we have described a relationship of the $G$-equivariant Conley index of the orbit $G(x_0)$ and $G_{x_0}$-equivariant Conley index of the $\{x_0\}$ considered as an isolated critical point of the  potential restricted to the orthogonal complement of 
$T_{x_0}G(x_0)$ at $x_0.$ Theorem \ref{smash theorem for sci} is an infinite-dimensional generalization of a combination of Theorems \ref{hsmadmco} and \ref{cio}. This theorem will play important role in the proof of the main result of our paper. In the last subsection of Section \ref{preliminaria} we have proved the splitting theorem, see Theorem \ref{splittinglemma}, which plays a crucial role in the study of isolated degenerate critical points.

Section \ref{results} is devoted to the proof of the main results of this article.  The study of periodic solutions of any period of system \eqref{newsys} is equivalent to the study of $2\pi$-periodic solutions of a family of systems, see \eqref{newsysp}. 

We have considered the solutions of system \eqref{newsysp}  as critical orbits of $G=(\Gamma \times \sone)$-invariant family of functionals $\Phi(q,\lambda)$ of class $C^2$ defined by formula \eqref{functional}. The necessary and sufficient conditions for the existence of local bifurcations of solutions of equation $\nabla_q \Phi(q,\lambda)=0$ have been proved in Section \ref{bifurcations}, see Theorems \ref{necessary}, \ref{sufficient}, respectively. In Section \ref{section} we study the $G_{q_0}$-equivariant Conley index  on the space orthogonal  to the orbit $G(q_0)$ at $q_0,$ see Lemmas  \ref{main-lemma0}, \ref{main-lemma}. In Subsection \ref{proof} we have proved the main result of  this paper. 

Section \ref{applications} contains the illustration of the abstract result of our article. We apply Theorem \ref{main-theo} to prove the existence of non-stationary periodic solutions of the Lennard-Jones and Schwarzschild problems, whose  potentials are $\Gamma=\sotwo$-invariant.

\section{Preliminary results}
\label{preliminaria}
\numsec 

In this section, for the convenience of the reader, we repeat  the relevant material from \cite{KBO,DIECK} without proofs, thus making our exposition self-contained. Moreover, we prove some preliminary results. Throughout this section $G$ stands for a compact Lie group.

\subsection{Groups and their representations }
\numsubsec
Denote  by $\subg$ the set of all closed subgroups of  $G.$ Two subgroups $H, H' \in \subg$ are said to be conjugate in $G$  if there is $g \in G$ such that $H=gH'g^{-1}.$ The conjugacy is an equivalence relation on $\subg.$ The class of $H \in \subg$ will be denoted by $(H)_G$ and the set of conjugacy classes will be denoted by $\subgc.$ Denote by  $\rho : G \to \on$  a continuous homomorphism. The space $\bR^n$
with the $G$-action defined by $G \times \bR^n \ni (g,x) \to \rho(g)x \in \bR^n$ is said to be a real,  orthogonal representation of $G$ which we write $\bV=(\bR^n,\rho).$  To simplify notation we write $gx$ instead of $\rho(g)x$ and $\bR^n$ instead of $\bV$ if the homomorphism is given in general.
 
If $x \in \bR^n$ then a group $G_x=\{g \in G : g x = x\} \in \subg$ is said to be the isotropy group of $x$ and  $G(x)=\{gx : g \in G\}$ is called the orbit through  $x$. It is known that the orbit $G(x)$ is a smooth $G$-manifold $G$-diffeomorphic to $G \slash G_x.$  An open subset $\Omega \subset \bR^n$ is said to be $G$-invariant if  $G(x) \subset \Omega$ for every $x \in \Omega$. 

Two orthogonal representations of $G,$ say $\bV=(\bR^n,\rho), \bV'=(\bR^n,\rho'),$ are equivalent (briefly $\bV \approx_G \bV'$) if there exists an equivariant linear isomorphism $L : \bV \to \bV'$ i.e. the isomorphism $L$ satisfying $L(gx)=gL(x)$ for any $g \in G, x \in \bR^n.$ 
Put $D(\bV)=\{x \in \bV : \| x \| \leq 1\},\, S(\bV)=\partial D(\bV)$, $S^{\bV}=D(\bV) \slash S(\bV)$ and $B_{r}(\bV)=\{x \in \bV : \| x \| < r\}$. Since the representation $\bV$ is orthogonal, these sets are $G$ invariant. 

Denote by $\bR[1,m],\, m \in \bN,$   a two-dimensional representation of $\sone=\{z \in \bC : \mid z\mid=1\}$ with an action of   $\sone$
given by $(\Phi(e^{i\phi}),(x,y)) \longrightarrow (\Phi(e^{i\phi}))^m(x,y)^T,$ where  $\Phi(e^{i\phi})=
\left[\begin{array}{lr}
\cos \phi & -\sin \phi\\
\sin \phi&\cos \phi
\end{array}\right].$
For $k, m \in \mathbb{N}$ we denote by $\mathbb{R}[k,m]$ the direct sum of $k$ copies of $\bR[1,m]$, we also
denote by $\mathbb{R}[k,0]$  the  $k$-dimensional trivial representation of $\sone.$ The following classical result gives a complete classification (up to an equivalence) of finite-di\-men\-sio\-nal representations of $\sone$, see \cite{ADM}.

\begin{Theorem}\label{soneRepr}
\label{tk} If  $\bV$ is a representation of $\sone$, then there are finite sequences $\{k_i\},\, \{m_i\}$ satisfying $m_i\in \{0\}\cup \mathbb{N},\,  k_i\in \mathbb{N}$ for $1\le i \le r,\,   m_1  < m_2 < \dots < m_r$
such that $\bV$ is equivalent to $\ds \bigoplus^r_{i=1} \bR[k_i ,m_i]$ i.e.   $\bV \approx_{\sone} \ds \bigoplus^r_{i=1} \bR[k_i,m_i]$. Moreover, the
equivalence class of $\bV$  is uniquely determined
by sequences $\{k_i\}, \{m_i\}$.
\end{Theorem}

Below we recall the notion of an admissible pair, which was introduced in \cite{PCHRYST}.

\bdf \label{admissible}
Fix $H \in \subg.$ A pair $(G,H)$ is said to be \textit{admissible} if for any $K_1,K_2 \in \subh$  the following condition is satisfied: 
$\text{ if } (K_1)_H \neq (K_2)_H \text{ then } (K_1)_G \neq (K_2)_G.$
\edf

Recall that if $\Gamma$ is a compact Lie group, then the pair $( \Gamma \times \sone,\{e\} \times \sone)$ is admissible, see Lemma 2.1 of \cite{PCHRYST}. This property will play a crucial role in the next section.

\subsection{$G$-equivariant maps}
\numsubsec
Let  $(\bV,\langle\cdot,\cdot\rangle)$ be an orthogonal representation of $G.$ Fix an open $G$-invariant subset $\Omega\subset\bV$.

\bdf
A map $\phi : \Omega \to \bR$ of class $C^k$ is called \textit{ $G$-invariant $C^k$-potential}, if $\phi(gx)=\phi(x)$ for every $g \in G$ and $x \in \Omega$. The set of $G$-invariant $C^k$-potentials will be denoted by $C^k_{G}(\Omega,\bR)$.
\edf

\bdf 
A map $\psi : \Omega \to \bV$ of the class $C^{k-1}$ is called \textit{$G$-equivariant $C^{k-1}$-map}, if $\psi(gx)=g \psi(x)$ for every $g \in G$ and  $x \in \Omega.$ The set of $G$-equivariant $C^{k-1}$-maps  will be denoted by $C^{k-1}_{G}(\Omega,\bV).$
\edf

 Fix $\vp \in C^2_G(\Omega,\bR)$ and denote by $\nabla \vp, \nabla^2 \vp$ the gradient and the Hessian of $\vp,$ respectively.
For $x_0 \in \Omega$ denote by  $\morse(\nabla^2 \vp(x_0))$ the Morse index of the  Hessian of $\vp$ at $x_0$  i.e. the sum of the  multiplicities of negative eigenvalues of the symmetric matrix $\nabla^2 \vp(x_0).$

\br \label{eqgrad} It is clear that if $\vp \in C^k_G(\Omega,\bR),$ then $\nabla \vp \in C^{k-1}_G(\Omega,\bV)$. Moreover,  if $x_0 \in (\nabla \vp)^{-1}(0),$ then $G(x_0)   \subset (\nabla \vp)^{-1}(0)$  i.e. the $G$-orbit of a critical point consists of critical points.  If $\nabla \vp(x_0)=0$ then  $\nabla \vp(\cdot)$ is fixed on $G(x_0).$ That is why $T_{x_0} G(x_0) \subset\ker \nabla^2 \vp(x_0) $ and consequently $\dim \ker \nabla^2 \vp(x_0) \geq \dim T_{x_0} G(x_0)=\dim G(x_0).$ 
\er

\subsection{Equivariant Conley index and equivariant Euler characteristic} \label{conleyIndexSection}
\numsubsec 
Denote by $\cF_{\ast}(G)$ the category of 
finite pointed $\gcw$-complexes, see \cite{DIECK}, where morphisms are continuous $G$-equivariant maps preserving base points.  By $\cF_{\ast}[G]$ we denote the set of
$G$-homotopy types of elements of $\cF_{\ast}(G),$ where $[\bX]_G\in \cF_{\ast}[G]$ (or $[\bX]$ when no confusion can arise) denotes a $G$-homotopy type of the pointed $G$-CW complex $\bX\in\cF_{\ast}(G)$. If $\bX$ is a $G$-CW-complex without a base point, then we denote by $\bX^+$ a pointed $G$-CW-complex $\bX^+=\bX\cup\{\ast\}$. A finite-dimensional $G$-equivariant Conley index of an isolated invariant set $S$ under a $G$-equivariant vector field $\vartheta$ will be denoted as $\cig(S,\vartheta)$, see \cite{BARTSCH,FLOER,GEBA,SMWA} for the definition. Recall that $\cig(S,\vartheta)\in\cF_{\ast}[G]$, see \cite{GEBA}. 

Below we present the infinite-dimensional extension of the equivariant Conley index due to Izydorek \cite{IZYDOREK} which requires the notion of equivariant spectra, see also \cite{GEIZPR,RYBAKOWSKI}.

Let $\xi=(\bV_n)_{n=0}^{\infty}$ be a sequence of finite-dimensional orthogonal representation of $G$.
\bdf\label{G-spectrum}
A pair $\cE(\xi)=\left((\cE_n)_{n=n(\cE(\xi))}^{\infty}, (\varepsilon_n)_{n=n(\cE(\xi))}^{\infty}\right)$, where $n(\cE(\xi))\in\bN$, is called a $G$-spectrum of type $\xi$ if
\begin{enumerate}
\item $\cE_n\in\cF_{\ast}(G)$ for $n\geq n(\cE(\xi))$,
\item $\varepsilon_n\in Mor_G(S^{\bV_n}\wedge \cE_n,\cE_{n+1})$ for $n\geq n(\cE(\xi))$,
\item there exists $n_1(\cE(\xi))\geq n(\cE(\xi))$ such that for $n \geq n_1(\cE(\xi))$, $\varepsilon_n$ is a $G$-homotopy equivalence.
\end{enumerate}
\edf
The set of $G$-spectra of type $\xi$ is denoted by $GS(\xi)$.

\bdf
A $G$-map of $G$-spectra $\cE(\xi),\cE '(\xi)$ is a sequence of maps \\
$f=(f_n)_{n=n_0}^{\infty}:\cE(\xi)\to\cE'(\xi)$, where $n_1\geq\max(n_1(\cE(\xi)),n_1(\cE '(\xi)))$, such that
\begin{enumerate}
\item $f_n\in Mor_G(\cE_n,\cE_n')$ for $n\geq n_1$,
\item $G$-maps $f_{n+1}\circ\varepsilon_n$ and $\varepsilon_n'\circ S^{\bV_n}f_n$ are $G$-homotopic for every $n\geq n_1$, where $S^{\bV}f_n$ denotes a suspension of $f_n$.
\end{enumerate}
\edf

 Two $G$-maps $f,g:\cE(\xi)\to\cE'(\xi)$ are $G$-homotopic if there exists $n_1\geq n_0$ such that $f_n,g_n:\cE\to\cE'$ are $G$-homotopic for $n\geq n_1$. Following this definition in a natural way we understand a $G$-homotopy equivalence of two spectra $\cE(\xi),\cE'(\xi)$. The $G$-homotopy type of a $G$-spectrum $\cE(\xi)$ will be denoted by $[\cE(\xi)]_G$ (or shorter $[\cE(\xi)]$) and the set of $G$-homotopy types of $G$-spectra by $[GS(\xi)]$ or simply $[GS]$ when $\xi$ is fixed or is not known yet.

\br
It follows from Definition \ref{G-spectrum} that the $G$-homotopy type $[\cE(\xi)]$ of spectrum $\cE(\xi)=\left((\cE_n)_{n=n(\cE(\xi))}^{\infty}, (\varepsilon_n)_{n=n(\cE(\xi))}^{\infty}\right)$ depends only on the sequence $(\cE_n)_{n=n_1(\cE(\xi))}^{\infty}$.
\er

Define an infinite-dimensional generalization of the equivariant Conley index. Since in this article we consider periodic solutions of second order systems as critical orbits of invariant functionals,  we remind the definition of the Conley index  only for gradient operators. Let $(\bH,\langle\cdot,\cdot\rangle)$ be an infinite-dimensional orthogonal Hilbert representation of  $G$. 

Let $L:\bH\to\bH$ be a linear, bounded, self-adjoint and $G$-equivariant operator with spectrum $\sigma(L)$ such that

\be \label{B-properties}
\item [(B.1)] $\bH=\overline{\bigoplus_{n=0}^{\infty} \bH_n}$, where all subspaces $\bH_n$ being mutually orthogonal representations of $G$ of finite dimension,
\item [(B.2)]$L(\bH_n)=\bH_n$ for all $n\geq 0$,
\item [(B.3)]$0$ is not an accumulation point of $\sigma (L)$.
\ee

Put $\bH^n:=\bigoplus_{k=0}^n\bH_k$ and denote by $P_n:\bH\to\bH^n$ the orthogonal projection onto $\bH^n$. Moreover denote by $\bH_k^+$ the subspace of $\bH_k$ corresponding to the positive part of spectrum of $L$.
Consider a functional $\Phi:\bH\to\bR$ such that $\nabla\Phi(x)=Lx+\nabla K(x)$ where $\nabla K\in C^1_H(\bH,\bH)$ is completely continuous. Denote by $\vartheta$ a \gls-flow, see Definition 2.1 of \cite{IZYDOREK}, generated by $\nabla\Phi$. Let be $\mathcal{O}$ an isolating $G$-neighborhood for $\vartheta$ and put $\cN=Inv_{\vartheta}\cO$. Set $\xi=(\bH_k^+)_{k=1}^{\infty}$. Let $\Phi_n:\bH^n\to\bR$ be given by $\Phi_n=\Phi_{\mid \bH^n}$ and $\vartheta_n$ denotes the $G$-flow generated by $\nabla\Phi_n$. Note that $\nabla\Phi_n(x)=Lx+P_n\circ\nabla K(x)$. Choose sufficiently large $n_0$ such that for $n\geq n_0$ the set $\cO_n:=\cO\cap\bH^n$ is an isolating $G$-neighborhood for the flow $\vartheta_n$. Then the set $Inv_{\vartheta_n}(\cO_n)$ admits a $G$-index pair $(Y_n,Z_n)$. 

We define a spectrum $\cE(\xi):=(Y_n/Z_n)_{n= n_0}^{\infty}$. 
Then the equivariant Conley index of $\cO$ with respect to the flow $\vartheta$ is given by
$
\scig(\cO,\vartheta):=[\cE(\xi)]\in [GS].
$
Sometimes we will write a vector field and isolated invariant set instead a flow and isolating neighborhood i.e. $\scig(\cN,\nabla\Phi)$.

Let $\left(U(G),+,\star)\right)$ be the Euler ring of  $G$, see 
\cite{DIECK} for the definition and properties of this ring. 
Let us briefly recall that the Euler ring $U(G)$ is commutative, generated by $\chig({G/H}^+)$, $(H)\in\subgc$ with the unit $\bI=\chig({G/G}^+)$, where $\chig:\cF_{\ast}[G]\to U(G)$ is the universal additive invariant for finite pointed $G$-CW-complexes known as the equivariant Euler characteristic.

\br\label{propOfChi}
Below we present some properties of the equivariant Euler characteristic $\chi_G(\cdot)$.
\begin{itemize}
\item  For  $\bX, \bY \in \cF_{\ast}(G)$ we have: $\chi_G(\bX)+\chi_G(\bY)=\chi_G(\bX\vee\bY)$ and $\chi_G(\bX)\star\chi_G(\bY)=\chi_G(\bX\wedge\bY)$
\iffalse 
\item If $\bX \in \mathcal{F}_{\ast}(G)$  then 
\beq \label{euler} \ds  \chi_G(\bX) = \sum_{(K)_G \in
\subgc} n^G_{(K)_G}(\bX) \cdot \chi_G\left(G\slash K^+\right), 
\eeq  where $\ds n^G_{(K)_G}(\bX)=\sum_{i=0}^{\infty} (-1)^{i} n(\bX,(K)_G,i)$ and $n(\bX,(K)_G,i)$ is the number of $i$-cells of type $(K)_G$ of $\bX$, see \cite{[DIECK]} for details.
\fi

\item If $\bW$ is an orthogonal representation of $G$ then $\chi_G(S^{\bW})$ is   invertible in $U(G)$, see \cite{DGR}.

\item If $\bW$ is a  representation of $\sone$ and $\bW \approx_{\sone} \bR[k_0,0] \oplus \bR[k_1,1] \oplus \ldots \bR[k_r,m_r]$, see Theorem \ref{soneRepr}, then
\beq \label{ches}
\chi_{\sone}\left(S^{\bW}\right)=\chi_{\sone}\left(S^{ \bR[k_0,0] \oplus \bR[k_1,1] \oplus \ldots \bR[k_r,m_r]}\right)=(-1)^{k_0}\left(\bI - \sum_{i=1}^r k_i  \chi_{\sone}({\sone \slash \bZ_{m_i}}^+)\right).
\eeq
The above equality has been proved in \cite{GLG}. See \cite{GLG} for more properties of the Euler ring $U(\sone)$ and the Euler characteristic $\chi_{\sone}$.

\end{itemize}
\er

 There is a natural extension of the equivariant Euler characteristic for finite pointed $G$-CW-complexes to the category of $G$-equivariant spectra due to Go{\l}\c{e}biewska and Rybicki \cite{GORY1}. 

 Let $\xi=(\bV_n)_{n=0}^{\infty}$ and put $\bV^n=\bV_0\oplus\bV_1\oplus\ldots\oplus\bV_n$, for $n\geq 0$. Recall that due to Remark \ref{propOfChi} an element $\chig(S^{\bV^n})$ is invertible in the Euler ring $U(G)$ and define a map $\upsg:[GS(\xi)]\to U(G)$ by the following formula
\beq \label{gchargspe}
\upsg([\cE(\xi)])=\lim_{n\to\infty} \left(\chig\left(S^{\bV^{n-1}}\right)^{-1}\star \chig(\cE_n)\right).
\eeq

\br\label{stabilization of upsilon}
It was shown in \cite{GORY1}  that $\upsg$ is well-defined. In fact 
\beq
\upsg([\cE(\xi)])=\chig\left(S^{\bV^{n_1(\cE)-1}}\right)^{-1}\star \chig(\cE_{n_1(\cE)})),
\eeq
where $n_1(\cE)=n_1(\cE(\xi))$ comes from Definition \ref{G-spectrum}.
\er

\br\label{Upsilon on finite complex}
Note that a finite pointed $G$-CW-complex $\bX$ can be considered as a \textit{constant} spectrum $\cE(\xi)$ where $\cE_n=\bX$ for all $n\geq 0$ and $\xi$ is a sequence of trivial, one-point representations. Then
$\ds
\upsg([\bX])=\upsg([\cE(\xi)])=\bI^{-1}\star\chig([\bX])=\chig([\bX]).
$
Therefore we can treat $\scig$ and $\upsg$ as natural extensions of $\cig$ and $\chig$ respectively.
\er

 By Theorems 3.1, 3.5  of \cite{GORY1} we obtain the following product formula.

\bt\label{product formula}
If $\cN_1,\,\cN_2$ are isolated $G$-invariant sets for the local $G$-$\mathcal{LS}$ flows generated by $\nabla\Psi_1$ and $\nabla\Psi_2$ respectively then
\[
\upsg\left(\scig\left(\cN_1\times\cN_2,(\nabla\Psi_1,\nabla\Psi_2)\right)\right)=\upsg\left(\scig\left(\cN_1,\nabla\Psi_1\right)\right) \star \upsg\left(\scig\left(\cN_2,\nabla\Psi_2\right)\right)
\]
\et

\subsection{How to distinguish two equivariant Conley indexes?}
\numsubsec
Let   $(\bV,\langle\cdot,\cdot\rangle)$ be a finite-dimen\-sio\-nal orthogonal representation of $G.$ Throughout this subsection   $\Omega\subset\bV$ stands for an open and $G$-invariant subset.

Fix a potential $\vp \in C^2_G(\Omega,\bR)$ and $x_0 \in (\nabla\vp)^{-1}(0)$. Suppose that $G(x_0)$ is an isolated orbit of critical points of $\vp$.
Our aim is to simplify the computation of the $G$-equivariant Conley index $\ci_G(G(x_0),-\nabla \vp) \in \cF_{\ast}[G]$ of the orbit  $G(x_0)$. We express it in terms of the $G_{x_0}$-equivariant Conley index $\ci_{G_{x_0}}(\{x_0\},-\nabla \phi) \in \cF_{\star}[G_{x_0}]$ of the critical point $x_0$ of the potential $\vp$ restricted to the space orthogonal to the orbit $G(x_0)$ i.e. $\phi=\vp_{\mid \tangent}$.

\br
Note that since $\bV$ is an orthogonal representation of $G$, the space $\tangent$ is an orthogonal representation of $G_{x_0}$.
\er

\bdf\label{D:twistedProduct}
Fix $H \in \subg$ and a $H$-space $\bY$. The product $G \times \bY$ carries $H$-action $(h,(g,y)) \to (gh^{-1},hy).$ The orbit space of $H$-action is denoted by $G \times_H \bY$ and called the twisted product over $H.$ We observe that $G \times_H \bY$ is a $G$-space with $G$-action defined by $(g',[g,x]) \to [g'g, x].$
\edf

\bdf\label{D:smashProduct} Let $\bY$ be a pointed $H$-space with a base point $\ast.$ Denote by  $G^+$  the group $G$ with  disjoint $G$-fixed base point $\ast$ added. Define the smash product of $G^+$ and $\bY$ by \linebreak $G^+ \wedge \bY=G^+ \times \bY \slash G^+ \vee \bY =G \times \bY \slash G \times \{\ast\}.$ The group $H$ acts on the pointed space $G^+ \wedge \bY$ by $(h,[g,y]) \to [gh^{-1},hy].$ The orbit space is denoted by $G^+ \wedge_H \bY$ and called the smash over $H,$ see \cite{DIECK}. A~formula  $(g',[g,y]) \to [g'g,y]$ induces $G$-action so that $G^+ \wedge_H \bY$ becomes a pointed $G$-space.
\edf

 Let $\bY\in\cF_{\ast}(H)$. The following theorem gives an interesting and very useful relation between Euler characteristics of $[\bY]$ and $[G^+ \wedge_H \bY]$. The proof of this theorem can be found in \cite{PCHRYST}.
 
\bt \label{hsmadmco}
Fix $H \in \subg$ such that the pair $(G,H)$ is admissible.  If $[X]_H, [Y]_H  \in \cF_{\ast}[H]$ and $\ds \chih([X]_H) \neq  \chih([Y]_H) \in U(H),$ then 
$\chig\left([G^+ \wedge_H X]_G\right) \neq \chig\left([G^+ \wedge_H Y]_G\right) \in U(G).$  In other words the map  $U(H) \ni \chih([X]_H) \to \chig\left([G^+ \wedge_H X]_G\right) \in  U(G)$ is injective.
\et

 In the theorem below we  express the $G$-equivariant Conley index  $\ci_G(G(x_0),-\nabla \vp)$ in  terms of the $H$-equivariant Conley index  $\ci_H  (\{x_0\},-\nabla \phi)$.

\bt \label{cio} 
Let $\vp \in C^2_G(\Omega,\bR)$. Suppose that $G(x_0)$ is an isolated orbit of critical points of $\vp$. Let $\phi \in C^2_{G_{x_0}}(\Omega \cap \tangent,\bR)$ be the restriction of $\vp$ to the space $\Omega \cap \tangent$. Then
$$\ci_G(G(x_0),-\nabla \vp)= G^+ \wedge_{G_{x_0}}\ci_{G_{x_0}}  (\{x_0\},-\nabla \phi) \in \cF_{\ast}[G].$$
\et
\begin{proof} To simplify notations we put $H=G_{x_0}.$
Firstly we express the $G$-index pair $(\cN,\cL)$ of the orbit $G(x_0)$ in terms of the twisted product over $H$ of the $H$-index pair $(N,L)$ of $x_0$. In fact $(\cN,\cL)=(G \times_H N,G \times_H L)$ is a $G$-index pair of the isolated invariant set $G(x_0)$, see \cite{KBO}. Consequently we obtain the  following equality 
\beq 
\label{gci} \cC\cI_G(G(x_0),-\nabla \vp)= ([(G \times_H N)\slash (G  \times_H L)]_G,\ast) \in \cF_{\ast}[G].
\eeq
It is clear that the spaces $(G \times N) \slash (G \times L)$ and $(G \times (N \slash L)) \slash (G \times \{\ast\})=G^+ \wedge (N \slash L)$ are homeomorphic. Moreover the $H$-action on all sets above is the same. Therefore the $H$-orbit spaces with given $G$-action  are $G$-homeomorphic i.e.
\[
(G\times_H N)\slash (G\times_H L)\approx_G G^+ \wedge_H (N\slash L)
\]
and as a consequence
\[
([(G \times_H N)\slash (G  \times_H L)]_G,\ast)=([G^+ \wedge_H (N \slash L)]_G,\ast)=G^+ \wedge_H ([N\slash L]_H,\ast),
\]
which completes the proof.
\end{proof}

Let $\bH=\overline{\bigoplus_{n=0}^{\infty} \bH_n}$ be a representation of $G.$ Consider two functionals $\vp_1,\vp_2\in C^2_G(\bH,\bR)$ such that $\nabla\vp_i=Lx+\nabla K_i(x)$, where $\nabla K_i\in C^1_G(\bH,\bH)$ is completely continuous, $i=1,2$ satisfying the conditions (B.1)--(B.3) described previously in Subsection \ref{conleyIndexSection}.

\bt\label{smash theorem for sci}
Let  $G(x_1), G(x_2)$ be isolated orbits of critical points of the potentials $\vp_1$ and $\vp_2$, respectively. Moreover, assume that $G_{x_1}=G_{x_2}(=H)$. If the pair $(G,H)$ is admissible and
$\upsh(\sci_H(\{x_1\},-\nabla \phi_1)) \neq \upsh(\sci_H(\{x_2\},-\nabla \phi_2)) \in U(H)$ where $\phi_i={\vp_i}_{\mid T_{x_i}^{\perp} G(x_i)}$ then 
\[\upsg(\sci_G(G(x_1),-\nabla \vp_1)) \neq  \upsg(\sci_G(G(x_2),-\nabla \vp_2)) \in U(G).
\]
\et
\begin{proof}
Since $\phi_1,\,\phi_2$ are in the form of compact perturbation of the same linear operator, the Conley indexes $\sci_H(\{x_1\},-\nabla \phi_1)$ and $\sci_H(\{x_2\},-\nabla \phi_2)$ are the homotopy types of spectra of the same type $\xi=(\bV_n)_{n=0}^{\infty}$. Denote these spectra by $\cE_1(\xi)$ and $\cE_2(\xi)$ respectively. Choose $n_1\geq \max\{n_1(\cE_1(\xi)),n_1(\cE_2(\xi))\}$, where $n_1(\cE_i(\xi)),\,i=1,2$ comes from Definition \ref{G-spectrum}. Therefore by Remark \ref{stabilization of upsilon} we obtain
\[
\upsh(\sci_H(\{x_i\},-\nabla \phi_i))=\chih\left(S^{\bV^{n-1}}\right)^{-1}\star \chih(\cih(\{x_i\},-\nabla\phi_i^{n}))
\]
for any $n\geq n_1$, where $\phi_i^n={\phi_i}_{\mid \bH^n}$, $\bH^n=\bigoplus_{i=0}^n \bH_i$ and $\bV^n=\bigoplus_{i=0}^n \bV_n$. Therefore the assumption $\upsh(\sci_H(\{x_1\},-\nabla \phi_1)) \neq \upsh(\sci_H(\{x_2\},-\nabla \phi_2))$ is equivalent to 
\beq\label{temp101}
\chih(\cih(\{x_1\},-\nabla\phi_1^{n}))\neq\chih(\cih(\{x_2\},-\nabla\phi_2^{n}))
\eeq
for any $n\geq n_1$.

The same reasoning can be performed for $\vp_1$ and $\vp_2$ obtaining similar formula for $n\geq n_2$ and $G$-equivariant Euler characteristics.

\noindent With all the above it is sufficient to show that  $\chih(\cih(\{x_1\},-\nabla\phi_1^{n}))\neq\chih(\cih(\{x_2\},-\nabla\phi_2^{n}))$ implies $\chig(\cig(G(x_1),-\nabla\vp_1^n))\neq\chig(\cig(G(x_2),-\nabla\vp_2^n))$, where $\vp_i^n={\vp_i}_{\mid\bH^n},\,i=1,2$ and $n\geq \max\{n_1,n_2\}$.

Fix $n\geq \max\{n_1,n_2\}$. Since $\phi_i^n={\vp_i^n}_{\mid T_{x_i}^{\perp} G(x_i)}$ we can apply Theorem \ref{cio} to obtain
\beq\label{temp100}
\ci_G(G(x_i),-\nabla \vp_i^n)= G^+ \wedge_H\ci_H  (\{x_i\},-\nabla \phi_i^n). 
\eeq
Combining  equalities \eqref{temp101}, \eqref{temp100} with  Theorem \ref{hsmadmco} we obtain
\[
\chig(\cig(G(x_1),-\nabla\vp_1^n))\neq\chig(\cig(G(x_2),-\nabla\vp_2^n)),
\]
which implies that 
$\upsg(\sci_G(G(x_1),-\nabla \vp_1)) \neq  \upsg(\sci_G(G(x_2),-\nabla \vp_2)) \in U(G).
$
\end{proof}

\subsection{Equivariant splitting lemma}
\label{splitting}
\numsubsec 
 
Let $H$ be a compact Lie group and let $(\bV,\langle\cdot,\cdot\rangle)$ be an orthogonal Hilbert representation  of $H$  with an invariant scalar product $\langle\cdot,\cdot\rangle$. Assume additionally   that $\dim \bV^H < \infty.$ Here and subsequently, $\Omega \subset \bV$ stands for an open and invariant subset of $\bV$ such that $0 \in \Omega.$

 Consider a functional $\Psi\in C^2_H(\Omega,\bR)$ given by the formula
\beq\label{psidef}
	\Psi(x)=\frac{1}{2}\langle Ax,x \rangle+\zeta(x),
\eeq
which satisfies the following assumptions
\be
\item [(F.1)] $A:\bV\to\bV$ is a $H$-equivariant self-adjoint linear Fredholm operator,
\item [(F.2)] $0\in\sigma(A)$,
\item [(F.3)] $\ker A \subset \bV^H$,
\item [(F.4)] $\nabla\zeta:\bV\to \bV$ is a $H$-equivariant, compact operator,
\item [(F.5)] $\nabla\zeta(0)=0$ and $||\nabla^2\zeta(x)||\to 0$ as $||x||\to 0$,
\item [(F.6)] $0 \in \Omega$ is an isolated critical point of $\Psi$.
\ee

Denote by $\ker A$ and $\im A$ the kernel and the image of $\nabla^2\Psi(0)=A,$ respectively. Notice that both, $\ker A$ and $\im A,$ are orthogonal representations of $H.$ Moreover, $\ker A$ is  finite dimensional and trivial representation of $H.$ Since $A$ is self-adjoint $\bV=\ker A \oplus \im A.$ Denote by $P:\bV\to \ker A$ and $Q=Id-P:\bV\to \im A$ the $H$-equivariant, orthogonal projections. 
Put $x=(u,v)$, where $u \in \ker A$ and $v \in \im A.$
The theorem below is an equivariant version of the implicit function theorem. A proof of the theorem below  can be found in \cite{FUGORY}.

\bt\label{implicitthm}
Let $\Psi\in C^2_H(\Omega,\bR)$ be a functional defined by \eqref{psidef} which satisfies assumptions $(F.1)$--$(F.5)$. Then there exists $\varepsilon_0>0$ and  $C^1$-mapping $w : B_{\varepsilon_0}(\ker A) \to \im A \cap \bV^H$ such that
\begin{enumerate}
\item $w(0)=0$, $Dw (0)=0$,
\item $Q\nabla\Psi (u,v)=0$ for $u\in D_{\varepsilon_0}(\ker A)$ iff $v=w(u)$.
\end{enumerate}
\et

In the following theorem (known as the \emph{splitting lemma}) we prove the existence of equivariant homotopy which allows us to study the product (splitted) flow $(\nabla\varphi(u),Av)$ where $u\in\ker A,\,v\in\im A$ instead of the general $\Psi(x)=\frac{1}{2}\langle Ax,x \rangle+\zeta(x)$. Note that $A$ is an isomorphism on $\im A$ and therefore the study the second flow is standard. Moreover, we are able to describe the first flow $\nabla\varphi(u),$ see Remark \ref{grad_varphi}.

\bt\label{splittinglemma} Suppose that a functional $\Psi\in C^2_H(\Omega,\bR)$ is defined by formula \eqref{psidef} and satisfies assumptions $(F.1)$--$(F.6)$. Then, there exists $\varepsilon_0>0$ and $H$-equivariant  gradient homotopy $\nabla \cH:(B_{\varepsilon_0}(\ker A) \times B_{\varepsilon_0}(\im A)) \times [0,1] \to \bV$ satisfying the following conditions:
\begin{enumerate}
\item $\nabla \cH((u,v),t)=Av-\nabla \xi_t(u,v)$, for $t\in [0,1]$, where $\nabla \xi_t=\nabla \xi(\cdot,t)$ and $\nabla \xi:\bV\times [0,1]\to\bV$ is  compact and  $H$--equivariant.
\item $(\nabla \cH)^{-1}(0)\cap (B_{\varepsilon_0}(\ker A)\times B_{\varepsilon_0}(\im A))\times[0,1]=\{0\}\times[0,1]$ i.e. $0$ is an isolated critical point of $\nabla\cH(\cdot,t)$ for any $t\in[0,1]$.
\item $\nabla \cH((u,v),0)=\nabla \Psi (u,v)$.
\item There exists an $H$--equivariant, gradient mapping $\nabla \varphi:B_{\varepsilon_0}(\ker A) \to \ker A$ such that $\nabla \cH((u,v),1)=(\nabla\varphi(u),Av),$ for all $(u,v) \in B_{\varepsilon_0}(\ker A) \times B_{\varepsilon_0}(\im A).$
\end{enumerate}
\et
\begin{proof}
First of all, we define a family of potentials $\cH :(B_{\varepsilon_0}(\ker A)\times B_{\varepsilon_0}(\im A))\times[0,1]\to \bR$ by
\begin{equation*}
		\cH((u,v),t)=\frac{1}{2}\langle Av,v\rangle +\frac{1}{2} t(2-t)\langle Aw(u),w(u)\rangle
		+ t\zeta (u,w(u))+(1-t)\zeta(u,v+tw(u)),
\end{equation*}
where $w$ is obtained  from Theorem \ref{implicitthm}. This family of functionals was introduced firstly by Dancer \cite{DANCER1}. Since we consider, contrary to Dancer, an infinite-dimensional and symmetric case, we check all the details in the proof of this theorem. Since $A, \zeta, w$ are $H$-equivariant and $\bV$  is an orthogonal representation of $H$, the functional $\cH(\cdot,t)$ is  $H$-invariant. Therefore $\nabla \cH$ is a gradient, $H$-equivariant homotopy, where $\nabla \cH$ denotes the gradient of $\cH$ with respect to the coordinate $x\in\bV$.
Observe that
\beq\label{gradker}
\begin{split}
P\nabla \cH((u,v),t)=&t(2-t)[Dw(u)]^T Aw(u)+tP\nabla\zeta(u,w(u))+t[Dw(u)]^T Q\nabla\zeta(u,w(u))+\\
		&+(1-t)P\nabla\zeta(u,v+tw(u))+t(1-t)[Dw(u)]^T Q\nabla\zeta(u,v+tw(u))
	\end{split}
\eeq
and that 
\beq\label{gradim}
	Q\nabla \cH((u,v),t)=Av+(1-t)Q\nabla\zeta(u,v+tw(u)).
\eeq 

Now we are ready to complete the proof of this theorem.

 (1) Put $\xi((u,v),t)=H((u,v),t)-\frac{1}{2}\langle Ax,x \rangle$. Taking into account formulas \eqref{gradker}, \eqref{gradim} we obtain
 \[
 \begin{split}
 -\nabla\xi(u,v,t)=&t(2-t)[Dw(u)]^T Aw(u)+tP\nabla\zeta(u,w(u))+t[Dw(u)]^T Q\nabla\zeta(u,w(u))+\\
		&+(1-t)P\nabla\zeta(u,v+tw(u))+t(1-t)[Dw(u)]^T Q\nabla\zeta(u,v+tw(u))+\\
&+(1-t)Q\nabla\zeta(u,v+tw(u))
 \end{split}
 \]
and therefore $\nabla \cH((u,v),t)=Av-\nabla \xi(u,v,t)$. The map $\nabla\xi$ is $H$-equivariant since $\nabla \cH$ and $A$ are $H$-equivariant.
To prove that $\nabla\xi$ is compact recall that $\nabla\zeta$ is compact, $\ker A \subset \bV^H$ is  finite dimensional and $w$ is defined on $\ker A$. Moreover,
\begin{itemize}
\item superposition of compact and continuous mappings is compact,
\item continuous, finite dimensional mapping is compact,
\item continuous mapping defined on a finite dimensional Banach space is compact.
\end{itemize}

(2) Since $\nabla^2 \cH((0,0),t)_{| \im A}=A_{|\im A}$ is an isomorphism, by the implicit function theorem there exists a unique solution $\tilde{w}(u,t)$ of $Q\nabla \cH((u,v),t)=0$ on  $B_{\varepsilon_1}(\ker A)\times B_{\varepsilon_1}(\im A)$ for small enough $\varepsilon_1$ independent of $t\in[0,1]$.
Let $\varepsilon_0<\frac{\min\{\varepsilon,\varepsilon_1\}}{2}$, where $\varepsilon$ is defined in Theorem \ref{implicitthm}. \\ Fix $t \in [0,1]$ and suppose that  $\nabla \cH((u,v),t)=0$ in $B_{\varepsilon_0}(\ker A)\times B_{\varepsilon_0}(\im A)$. Putting $(1-t)w(u)$ to \eqref{gradim} instead of $v$ we get
\[
Q \nabla \cH((u,(1-t)w(u)),t)=(1-t)\left(Aw(u)+Q\nabla\zeta(u,w(u))\right)=(1-t)Q\nabla \Psi(u,w(u))=0.
\]
By the uniqueness $\tilde{w}(u,t)=(1-t)w(u)$ is the only small solution of $Q\nabla \cH ((u,v),t)=0$ for each $t\in[0,1]$. 
Next putting $\tilde{w}(u,t)$ to \eqref{gradker} instead of $v$  we obtain
\begin{multline*}
P\nabla \cH ((u,(1-t)w(u)),t)=\\=t(2-t)[Dw(u)]^T Aw(u) + t P \nabla\zeta(u,w(u)) + t[Dw(u)]^T Q\nabla\zeta(u,w(u)) +\\
+(1-t) P\nabla\zeta(u,w(u))  +t(1-t)[Dw(u)]^T Q\nabla\zeta(u,w(u))=\\
=P\nabla\zeta(u,w(u))+t(2-t)[Dw(u)]^T \left(Aw(u) + Q\nabla\zeta(u,w(u))\right)=\\=P\nabla\zeta(u,w(u)) +t(2-t)[Dw(u)]^T Q\nabla\Psi(u,w(u) )= P\nabla\zeta(u,w(u)).
\end{multline*}
To summarize, if $\nabla \cH((u,v),t)=0$ (i.e. $Q\nabla \cH((u,v),t)=0$ and $P\nabla \cH((u,v),t)=0$) then $v=\tilde{w}(u,t)=(1-t)w(u)$ and $P\nabla\zeta(u,w(u))=0$. As a consequence we obtain
$$
\nabla\Psi(u,w(u))=Aw(u)+Q\nabla\zeta(u,w(u))+P\nabla\zeta(u,w(u))=\frac{1}{1-t}Q\nabla\cH((u,\tilde{w}),t)=0.
$$
Since $(0,0) \in \ker A \oplus \im A$ is an isolated critical point of $\Psi,$  $u=0$ and $v=\tilde{w}(0,t)=(1-t)w(0)=0$. Hence $(0,0)$ is an isolated critical point of  $\cH(\cdot,t)$ for any $t \in [0,1].$ Which completes the proof of (2).

To finish the proof of Theorem \ref{splittinglemma} it is enough to notice that:

 (3) $\nabla \cH((u,v),0)=\nabla\Psi(u,v)$, by the definition.

 (4) $\nabla \cH((u,v),1)=  P\nabla\zeta(u,w(u)) + Av= \left(\nabla\psi(u),Av\right),$ where $\nabla\psi(u)=P\nabla\Psi(u,w(u))=P\nabla\zeta(u,w(u))$.
\end{proof}

\br\label{grad_varphi}
The potential $\varphi: B_{\varepsilon_0}(\ker A) \to \bR$ is given by $\varphi(u)=\Psi(u,w(u))$. Indeed, for $z \in \ker A$, since
$
D\varphi(u)(z)=(D_u\Psi(u,w(u))+D_v\Psi(u,w(u)) Dw(u)) z
$
we obtain $$\langle\nabla \varphi(u),z\rangle=\langle P \nabla \Psi(u,w(u)),z)\rangle+\langle [Dw(u)]^T Q \nabla \Psi(u,w(u)),z\rangle = 
\langle P \nabla \Psi(u,w(u)),z\rangle.$$ 
\er

\section{main result}
\numsec \label{results}
In this section, using the equivariant bifurcation theory techniques, we prove the main result of this article i.e. the Symmetric Liapunov center theorem for minimal orbit, see Theorem \ref{main-theo}. 

We consider $\bR^n$ as an orthogonal representation of a compact  Lie group $\Gamma$. Denote by $\Omega \subset \bR^n$ an open and $\Gamma$-invariant subset. Fix $U \in C^2_{\Gamma}(\Omega,\bR)$ and  $q_0 \in \Omega$ a minimum of the potential $U$ such that isotropy group $\Gamma_{q_0}$ is trivial. Since $U$ is $\Gamma$-invariant, the orbit $\Gamma(q_0)$ consists of minima of $U$. Obviously $q_0$ is a critical point of $U$ and therefore $\Gamma(q_0)\subset (\nabla U)^{-1}(0)$ i.e. the orbit $\Gamma(q_0)$ consists of critical points of  $U$.

\br
Since $\Gamma_{q_0}$ is trivial, the orbit $\Gamma(q_0)$ is $\Gamma$-homeomorphic to $\Gamma \slash \Gamma_{q_0}=\Gamma$. For this reason it can happen that elements of this orbit are not isolated. For example, if $\Gamma=\sotwo$ acts freely on $\Omega$ then, the orbit $\Gamma(q_0)$ is $\sotwo$-homeomorphic to $\Gamma \slash \Gamma_{q_0}=\Gamma \slash \{e\} =\sotwo \approx  \sone$. Hence we can not treat $q\in\Gamma(q_0)$ as an isolated critical point of $U$. This is the reason to apply equivariant Conley index theory.
\er

 Note that the study of periodic solutions of any period of system \eqref{newsys} is equivalent to the study of $2\pi$-periodic solutions of the following  system 
\beq \label{sys} 
\left\{
\ba{rcl}  \ddot q(t) & = & -\lambda^2 \nabla U(q(t)), \\
q(0) & = & q(2\pi), \\
\dot q(0) & = & \dot q(2\pi).
\ea
\right.
\eeq 
The $2 \pi \lambda$-periodic solution of  \eqref{newsys} corresponds to $2\pi$-periodic solutions of \eqref{sys}. Since   $\Gamma(q_0) \subset (\nabla U)^{-1}(0),$ for every $\lambda > 0$ the orbit $\Gamma(q_0)$ consists of stationary solutions of   \eqref{sys}. 

\subsection{Variational setting}
\numsubsec 
In this article we treat solutions of  \eqref{sys} as  critical points of invariant functionals. This fact allows us to use equivariant bifurcation theory in order to prove our main result. Therefore we present variational setting for  family \eqref{sys}. 

Define 
\[
\h^1_{2\pi} = \{u : [0,2\pi] \rightarrow \bR^n : \text{ u is abs. continuous map, } u(0)=u(2\pi), \dot u \in L^2([0,2\pi],\bR^n)\}
\]
and a scalar product
\[
\ds \langle u,v\rangle_{\h^1_{2\pi}} = \int_0^{2\pi} (\dot u(t), \dot v(t)) + (u(t),v(t)) \; dt,
\]
where $(\cdot,\cdot)$ and $\| \cdot \|$ are the usual scalar product and norm in  $\bR^n,$ respectively. It is well known that $\left(\h^1_{2\pi},\langle\cdot,\cdot\rangle_{\h^1_{2\pi}}\right)$ is a separable Hilbert space. Moreover, it can be considered as an orthogonal representation of $G=\Gamma\times \sone$ where the action is given by
\[
G \times \h^1_{2\pi} \ni ((\gamma,e^{i \theta}),q(t)) \to \gamma q(t+\theta) \text{ mod } 2 \pi.
\]

It is known that solutions of system \eqref{sys} are in one to one correspondence with $\sone$-orbits of critical points of $\sone$-invariant potential $\Phi : \h^1_{2\pi}\times (0,\infty)  \to \bR$ of  class $C^2$ defined by
\beq  \label{functional}
\Phi(q,\lambda) = \int_0^{2\pi} \left( \frac{1}{2} \| \dot q(t) \|^2 -  \lambda^2U(q(t)) \right) \; dt,
\eeq
where $\lambda$ is considered as a parameter, see \cite{MAWI}. As $\bR^n$ is an orthogonal representation of $\Gamma$ and $U$ is $\Gamma$-invariant, the potential $\Phi$ is also $\Gamma$-invariant. Therefore $2\pi$-periodic solutions of system \eqref{sys} can be considered as critical orbits of  $G=(\Gamma \times \sone)$-invariant potential $\Phi$ i.e. as  solutions of the system 
\beq \label{gradf}
\nabla_q \Phi(q,\lambda)=0.
\eeq
Let $ \{ e_1,\ldots,e_n \} \subset \bR^n$ be the standard basis in $\bR^n.$ Define $\h_0=\bR^n, \h_k= \lin \{e_i \cos kt, e_i \sin kt: i=1,\ldots,n\}$ and note that 
\beq \label{space} 
\bH^1_{2\pi} = \overline{\h_0 \oplus \bigoplus_{k=1}^{\infty} \h_k} 
\eeq
and that the finite-dimensional spaces $\h_k,\,k=0,1,\ldots$ are orthogonal representations of $G$.

Note that the gradient $\nabla \Phi :
\h^1_{2\pi}\times(0,\infty)  \rightarrow \h^1_{2\pi}$ is a $G$-equivariant $C^1$-operator in the form of a compact perturbation of the
identity,  see \cite{PCHRYST} for more details. Summarizing, we will study the existence of $G$-orbits of critical points of $\Phi$. Let us underline that $\Phi$ satisfies assumptions (B.1)--(B.3) of Subsection \ref{conleyIndexSection}.

\br Assume that  $q_0 \in (\nabla U)^{-1}(0)$ and consider the linearization of the system \eqref{sys} at $q_0$ of the form 
% \beq% \label{sysl}  
\beq\label{linearizedSystem}
\left\{
\ba{rcl}  \ddot q(t) & = & - \lambda^2\nabla^2 U(q_0)(q-q_0), \\
q(0) & = & q(2\pi), \\
\dot q(0) & = & \dot q(2\pi).
\ea
\right. 
\eeq
The corresponding functional  $\Psi : \h^1_{2\pi} \rightarrow \bR$  is defined as follows
\beq\label{fun}
\Psi(q,\lambda)=\frac{1}{2} \|q\|^2_{\h^1_{2\pi}} + \langle \lambda^2\nabla^2 U(q_0)q_0,q \rangle_{\h^1_{2\pi}} - \frac{1}{2} \langle Lq,q\rangle_{\h^1_{2\pi}},
\eeq
where   $L : \h^1_{2\pi} \rightarrow \h^1_{2\pi}$ is a linear, self-adjoint, $G$-equivariant and compact operator, see \cite{PCHRYST} for details.
It is clear that $\nabla_q \Psi(q,\lambda)=q-Lq +\lambda^2\nabla^2 U(q_0)q_0.$

 If $q \in \h^1_{2\pi}$ is given by the Fourier series
$\ds q(t)= a_0   + \sum_{k =1}^{\infty} a_k \cdot   \cos   k  t   +
 b_k \cdot   \sin   k   t,$  we know that
\beq\label{fourier}
\nabla_q \Psi(q,\lambda)=  -\lambda^2\nabla^2 U(q_0)( a_0-q_0)   + \sum_{k=1}^{\infty}
 (Q(k,\lambda) \cdot a_k)   \cdot \cos  k t +  (Q(k,\lambda) \cdot b_k)  \cdot \sin k t,\eeq
 where $\ds Q(k,\lambda)=\left(\frac{k^2}{ k^2  + 1} Id - \frac{\lambda^2}{k^2 + 1} \nabla^2 U(q_0)\right),$ see Lemma 5.1.1 of \cite{FUGORY} for details.
\er

\subsection{Bifurcation from the critical orbit}
\numsubsec
\label{bifurcations} 
In this subsection we present  necessary and sufficient condition for the existence of local bifurcation of non-stationary periodic solutions from the family of stationary $\Gamma$-orbits.

Since $q_0 \in \h^1_{2\pi}$ is a constant function, $G(q_0)=\Gamma(q_0) \subset \bH_0 = \bR^n \subset \h^1_{2\pi}$ solves equation \eqref{gradf} for any $\lambda > 0.$ 
The set of solutions $\cT=G(q_0) \times (0,+\infty) \subset \h^1_{2\pi} \times (0,+\infty)$ of equation \eqref{gradf} is called a  family of trivial solutions of equation \eqref{gradf} while the set $\cN=\{(q,\lambda)\in \h^1_{2\pi} \times (0,+\infty)\setminus \cT\,:\, \nabla_q\Phi(q,\lambda)=0\}$ is called a family of non-trivial solutions.

\bdf
We say that the orbit $G(q_0)\times\{\lambda_0\}\subset\cT$ is an orbit of local bifurcation if $(q_0,\lambda_0)\in\overline{\cN}$, i.e. $(q_0,\lambda_0)$ is an accumulation point of nontrivial solutions of equation \eqref{gradf}.
\edf

 We present below the necessary condition for the existence of local bifurcation of   solutions of equation 
\beq\label{form1}
\nabla_q \Phi(q,\lambda)=0
\eeq
from the critical orbit $G(q_0)\times\{\lambda_0\}$. 
We define  $\Lambda=\{k\slash \beta: k \in \bN \text{ and } \beta^2 \in \sigma(\nabla^2 U(q_0)) \cap (0,+\infty)\}$.

\bt \label{necessary} (Necessary condition).
If   $G(q_0)\times\{\lambda_0\}$ is an orbit of local bifurcation of solutions of equation \eqref{form1} then $\ker \nabla^2_q\Phi(q_0,\lambda_0) \cap \overline{\bigoplus_{k=1}^{\infty} \h_k}\neq\emptyset$ i.e. $\lambda_0 \in \Lambda.$
\et
\begin{proof}
Suppose, contrary to our claim, that $\ker \nabla^2_q\Phi(q_0,\lambda_0) \cap \overline{\bigoplus_{k=1}^{\infty} \h_k} = \emptyset.$
Since $\nabla^2_q\Phi(q_0,\lambda_0)$ is  self-adjoint, $\h^1_{2\pi}=\ker \nabla^2_q\Phi(q_0,\lambda_0) \oplus \im \nabla^2_q\Phi(q_0,\lambda_0).$ Denote $q=(q_1,q_2)\in \ker \nabla^2_q\Phi(q_0,\lambda_0)\oplus \im \nabla^2_q\Phi(q_0,\lambda_0).$ Equation \eqref{form1} is equivalent to the following system
\beq\label{form2.1}
\nabla_{q_1}\Phi(q_1,q_2,\lambda)=0,
\eeq
\beq\label{form2.2}
 \nabla_{q_2}\Phi(q_1,q_2,\lambda)=0.
\eeq
Since $\nabla^2_{q^2}\Phi(q_0,\lambda_0)_{\mid  \im \nabla^2_q\Phi(q_0,\lambda_0)}$ is an isomorphism, applying Theorem  \ref{implicitthm} we obtain that \linebreak $(q_1,q_2,\lambda)=(q_1,q_2(q_1,\lambda),\lambda)$ is the only solution of \eqref{form2.2} in the neighborhood of $(q_0,\lambda_0)$, where $q_2: D_{\varepsilon} (\ker \nabla^2_q\Phi(q_0,\lambda_0)) \times [\lambda_0-\varepsilon, \lambda_0 + \varepsilon] \to \im  \nabla^2_q\Phi(q_0,\lambda_0)$ is a $G$-equivariant map. Therefore the study of equation \eqref{form1} is equivalent to study the equation
\beq\label{form3}
\nabla_{q_1}\Phi(q_1,q_2(q_1,\lambda),\lambda)=0.
\eeq
Now we are able to control the isotropy group of solutions of equation \eqref{form1} in the neighborhood of the orbit $G(q_0)\times\{\lambda_0\}$. Indeed, suppose that $(\widehat{q}_1,q_2(\widehat{q}_1,\widehat{\lambda}),\widehat{\lambda})$ is a solution of \eqref{form1}. Since $q_2$ is $G$-equivariant we obtain that $G_{\widehat{q}_1}=G_{(\widehat{q}_1,\widehat{\lambda})} \subset G_{q_2(\widehat{q}_1,\widehat{\lambda})}$ and therefore
\[
G_{(\widehat{q}_1,q_2(\widehat{q}_1,\widehat{\lambda}),\widehat{\lambda})}= G_{\widehat{q}_1} \cap G_{q_2(\widehat{q}_1,\widehat{\lambda})} \cap G_{\widehat{\lambda}}=G_{\widehat{q}_1},
\]
i.e. the isotropy groups of bifurcating solutions must coincide with isotropy groups of  elements of $\ker \nabla^2_q\Phi(q_0,\lambda_0)$. Since $G(q_0)=\Gamma(q_0)\subset \bH_0$ is an isolated orbit of constant solutions, the bifurcation can not occur in the direction of $\bH_0$. Finally, since  $\ker \nabla^2_q\Phi(q_0,\lambda_0)\subset \bH_0$ then $G(q_0)\times\{\lambda_0\}$ is not an orbit of local bifurcation, a contradiction.

 To complete the proof it is enough to show that
 $\ker \nabla^2_q\Phi(q_0,\lambda_0) \cap \overline{\bigoplus_{k=1}^{\infty} \h_k}\neq\emptyset$ if and only if $\lambda_0\in\Lambda$. The study of $\ker \nabla^2_q\Phi(q_0,\lambda_0)$ is equivalent to the study of the linearized system \eqref{linearizedSystem} and further, it is equivalent to the equation $\nabla_q \Psi(q,\lambda)=0$. By the equality \eqref{fourier}, the last equation has solutions in $\overline{\bigoplus_{k=1}^{\infty} \h_k}$ if and only if a matrix $\ds Q(k,\lambda)=\left(\frac{k^2}{ k^2  + 1} Id - \frac{\lambda^2}{k^2 + 1} \nabla^2 U(q_0)\right)$ is degenerate for some $k\in\bN$ i.e. $k=\lambda\beta$ for some $\beta^2\in\sigma(\nabla^2 U(q_0))\cap (0,\infty)$.
\end{proof}

The  theorem below provides the sufficient condition for the existence of local bifurcation in the terms of equivariant Conley index. This is a direct consequence of continuation property and homotopy invariance of equivariant Conley index, see \cite{IZYDOREK}.

\bt \label{sufficient} (Sufficient condition).
Under the assumptions above, if there exist $\lambda_1, \lambda_2$ such that $[\lambda_1,\lambda_2] \cap \Lambda = \{\lambda_0\}$   and
\beq
\sci_{G}\left(G(q_0),-\nabla\Phi(\cdot,\lambda_1)\right)\neq \sci_{G}\left(G(q_0),-\nabla\Phi(\cdot,\lambda_2)\right),
\eeq
then  $G(q_0)\times\{\lambda_0\}$ is an orbit of local bifurcation.
\et

\subsection{Equivariant Conley index on the orthogonal section} 
\numsubsec
\label{section}
In order to prove Theorem \ref{main-theo} we will study the existence of local bifurcation of solutions of equation \eqref{gradf} from the trivial family $\cT$. Additionally, we will control the minimal period of bifurcating solutions. The existence of local bifurcation will be a consequence of the change of the infinite-dimensional equivariant Conley index.

 Fix $\beta_{j_0}$ satisfying the assumptions of Theorem \ref{main-theo}, choose  $\varepsilon > 0$ sufficiently small and define $\ds \lambda_{\pm}=\frac{1\pm \varepsilon}{\beta_{j_0}}.$ Without loss of generality one can assume that $[\lambda_-,\lambda_+] \cap \Lambda = \{1/\beta_{j_0}\}.$ 
Since bifurcation does not occur at the level $\lambda_{\pm},$ the orbit $G(q_0)$ is isolated in $\left(\nabla \Phi(\cdot,\lambda_{\pm})\right)^{-1}(0)$. Therefore $G(q_0)$ is an isolated invariant set in the sense of the $G$-equivariant Conley index theory i.e. $\scig(G(q_0),-\nabla \Phi(\cdot,\lambda_{\pm}))$ is well-defined.

Let $\h \subset \h^1_{2\pi}$ be a linear subspace orthogonal  to $G(q_0)$ at $q_0$ i.e. $\h=T_{q_0}^{\perp} G(q_0) \subset \h^1_{2\pi}.$ Since  $G_{q_0}=\{e\} \times \sone,$ $\bH$ is an orthogonal representation of $\sone.$ Define a $\sone$-invariant functional of class  $C^2$ by $\Psi_{\lambda_{\pm}}=\Phi(\cdot,\lambda_{\pm})_{\mid \bH} : \bH \to \bR.$ Since $G(q_0) \subset \bH^1_{2\pi}$  is an isolated critical orbit of the $G$-invariant functional $\Phi(\cdot,\lambda_{\pm}),$  $q_0 \in \bH$ is an isolated critical point of $\sone$-invariant potential $\Psi_{\lambda_{\pm}}.$  Hence $q_0$ is an isolated invariant set in the sense of the $\sone$-equivariant Conley index theory defined in \cite{IZYDOREK} i.e. $\sci_{\sone}(\{q_0\},-\nabla \Psi_{\lambda_{\pm}})$ is defined. Let us underline that since we know the form of $\Phi$, $\Psi_{\lambda_{\pm}}$ satisfies the assumptions (F.1)--(F.6) given in Subsection \ref{splitting}.

Note that 
$$\bH_0=\bR^n= T^{\perp}_{q_0}\Gamma(q_0) \oplus  T_{q_0}\Gamma(q_0) \quad {\rm and} \quad  \ds \bH= T_{q_0}^{\perp} \Gamma(q_0) \oplus \overline{\bigoplus_{k=1}^{\infty} \bH_k} \subset \bH^1_{2\pi}.$$

 Let $\widetilde{U} : T_{q_0}^{\perp}\Gamma(q_0)\to\bR$ be given by the formula $\widetilde{U}(q)=U(q+q_0)$ and similarly $\widetilde{\Psi}_{\lambda_{\pm}}(q)=\Psi_{\lambda_{\pm}}(q+q_0)$. Since  $\nabla^2 \widetilde{\Psi}_{\lambda_{\pm}}(0)$ is  self-adjoint, $\bH=\ker \nabla^2 \widetilde{\Psi}_{\lambda_{\pm}}(0) \oplus \im  \nabla^2 \widetilde{\Psi}_{\lambda_{\pm}}(0).$ Moreover, by equality \eqref{fourier} $\ker \nabla^2\tpsil(0) =\ker \nabla^2 \widetilde{U}(0)\subset \bH_0$ is independent of $\lambda_{\pm}$.  As a consequence we obtain that  $\im \nabla^2 \widetilde{\Psi}_{\lambda_{\pm}}(0)=\im \nabla^2\Psi_{\lambda_{\pm}}(q_0)$ is independent of $\lambda_{\pm}.$ 
Because $\ker \nabla^2\widetilde{\Psi}_{\lambda_{\pm}}(0)$ and $ \im  \nabla^2\widetilde{\Psi}_{\lambda_{\pm}}(0)$ do not depend on $\lambda_{\pm}$ for abbreviation of notations  we write $\cN$ instead of $\ker \nabla^2\widetilde{\Psi}_{\lambda_{\pm}}(0)$ and $\cR$ instead of $\im  \nabla^2\widetilde{\Psi}_{\lambda_{\pm}}(0).$

Put $P:\bH \to\cN$ and $Q=Id-P : \bH \to \cR$ for the $H$-equivariant, orthogonal projections.
Note that the $\sone$-invariant  potential $\Pi_{\lpm} : \cR \to \bR$ of the linear vector field $\nabla^2{\Psi_{\lpm}}_{\mid \cR}(q_0)(q-q_0)$ is defined by $\ds \Pi_{\lambda_{\pm}}(q)=\frac{1}{2} \langle \nabla^2{\Psi_{\lpm}}_{\mid \cR}(q_0)(q-q_0),q-q_0\rangle_{\h^1_{2\pi}}.$

 In the following lemma we  reduce the computation of the $\sone$-equivariant Conley indexes of nonlinear maps to the linear case.

\bl\label{main-lemma0}
Under the above assumptions  $$\upssone\left(\sci_{\sone}(\{q_0\},-\nabla\Psi_{\lambda_+}\right)\neq \upssone\left(\sci_{\sone}(\{q_0\},-\nabla\Psi_{\lambda_-})\right)$$ 
$$ \textrm{if and only if}$$ $$\upssone\left(\sci_{\sone}\left(\{q_0\},-\nabla\Pi_{\lambda_+}\right)\right)\neq \upssone\left(\sci_{\sone}\left(\{q_0\},-\nabla\Pi_{\lambda_-}\right)\right).$$
\el
\begin{proof}
Since $\sci_{\sone}(\{q_0\},-\nabla \Psi_{\lambda_{\pm}})= \sci_{\sone}(\{0\},-\nabla \widetilde{\Psi}_{\lambda_{\pm}}),$  we will study $\sci_{\sone}(\{0\},-\nabla \widetilde{\Psi}_{\lambda_{\pm}})$ instead of $\sci_{\sone}(\{q_0\},-\nabla \Psi_{\lambda_{\pm}})$. 
Since the operator  $\widetilde{\Psi}_{\lambda_{\pm}}$ is defined by formula \eqref{psidef} and satisfies assumptions $(F.1)$-$(F.6)$ of Subsection \ref{splitting}, one  can combine Theorems \ref{product formula}, \ref{splittinglemma}  to obtain
\[
\upssone\left(\sci_{\sone}(\{0\},-\nabla \widetilde{\Psi}_{\lambda_{\pm}})\right)
=\upssone\left(\sci_{\sone}(\{0\},-\nabla\varphi_{\lambda_{\pm}})\right)\star \upssone\left(\sci_{\sone}(\{0\},-\nabla^2 \widetilde{\Psi}_{\lambda_{\pm}}(0)_{\mid \cR})\right),
\]
where $0=(0,0)\in \cN\oplus \cR $, $\varphi_{\lambda_{\pm}} : B_{\varepsilon_0}(\cN) \to \bR$, $\varphi_{\lambda_{\pm}}(u)=\widetilde{\Psi}_{\lambda_{\pm}}(u,w(u))$ and $\nabla\varphi_{\lambda_{\pm}}(u)=P\nabla \widetilde{\Psi}_{\lambda_{\pm}} (u,w(u))$ is $\sone$-equivariant. 

Now, since $\cN   \subset \bH_0=\bH^{\sone},$ applying Theorem \ref{implicitthm}  we obtain $(u,w(u))\in \bH^{\sone}=\bH_0$. Following the definition of $\widetilde{\Psi}_{\lambda_{\pm}}$ one can see that 
\beq \label{pote}
\varphi_{\lambda_{\pm}}(u)=\widetilde{\Psi}_{\lambda_{\pm}}(u,w(u))=-2\pi\lambda_{\pm}^2 \widetilde{U}(u,w(u)).
\eeq
Since $\{0\}$ is isolated in $(\nabla \varphi_{\lambda_{\pm}})^{-1}(0),$ one can choose   $0<\varepsilon<\varepsilon_0$ (where $\epsilon_0$ comes from Theorem \ref{implicitthm}) such that  $(\nabla \varphi_{\lambda_{\pm}})^{-1}(0)\cap D_{\varepsilon}(\cN)=\{0\}$. Recall that $v=w(u)$ is the only solution of $Q\nabla \widetilde{\Psi}_{\lambda_\pm}(u,v)=0$ for $u \in D_{\varepsilon}(\cN),$ see Theorem \ref{implicitthm}. 

Moreover $(0,0)=(0,w(0))\in \cN \oplus \left(\cR \cap \bH_0\right) \subset \bH_0$ is  an isolated critical point (minimum) of $\widetilde{U}$. Hence $0\in \cN$ is an isolated local maximum of $\varphi_{\lambda_{\pm}}$, see \eqref{pote}. Indeed, $u=0$ is an isolated critical point because  $\nabla\varphi_{\lambda_{\pm}}(u)=0$ implies $P(\nabla \widetilde{\Psi}_{\lambda_{\pm}})(u,w(u))=0$ and further $\nabla \widetilde{\Psi}_{\lambda_{\pm}}(u,w(u))=0$. However $(0,0)$ is the only solution of $\nabla\tpsil(u, w(u))=0$ on $D_{\varepsilon}(\cN)$. 

Since  $\cN \subset \bH^{\sone}$ is finite-dimensional, applying Remark \ref{Upsilon on finite complex} we obtain 
\[
\sci_{\sone}(\{0\},-\nabla\varphi_{\lambda})= \ci_{\sone}(\{0\},-\nabla\varphi_{\lambda})=\ci(\{0\},-\nabla\varphi_{\lambda}).
\]

Therefore 
$
\upssone\left(\sci_{\sone}(\{0\},-\nabla\varphi_{\lambda})\right)=\chi_{\sone}\left(\ci_{\sone}(\{0\},-\nabla\varphi_{\lambda})\right)=\chi\left(\ci(\{0\},-\nabla\varphi_{\lambda})\right)\cdot \bI \in U(\sone).
$
i.e.
the $\sone$-equivariant Euler characteristic of $\sci_{\sone}(\{0\},-\nabla\varphi_{\lambda})$ is generated by the identity in the Euler ring $U(\sone).$

 There exists a simple relation between the Euler characteristic of the Conley index $\chi(\ci(S,-\eta))$ of an isolated $\eta$-invariant set $S$ with an isolating neighborhood $N$ and the Brouwer degree $\deg (\nu,N)$, where $\eta$ is a local flow generated by the equation $\dot{x}=-\nu(x)$. \\ In fact $\chi(\ci(S,-\nu))=\deg (\nu,N),$ see  \cite{MCCORD,SRZEDNICKI} for details.

Rabinowitz proved  in \cite{RAB4} that the Brouwer index of an isolated critical point which is a minimum is equal to $1,$ see also \cite{AMANN}. This implies that the Brouwer index of an isolated maximum equals $(-1)^k$, where $k$ is the dimension of a space.

Hence,
\begin{equation*}
\begin{split}
\upssone\left(\sci_{\sone}(\{0\},-\nabla\varphi_{\lambda_{\pm}})\right)&= \chi\left(\ci(\{0\},-\nabla\varphi_{\lambda_{\pm}})\right)\cdot \bI=\\
&=\deg (\nabla\varphi_{\lambda_{\pm}}, D_{\varepsilon}(\cN))\cdot \bI=(-1)^{\dim \cN}\cdot \bI \in U(\sone).
\end{split}
\end{equation*}
As a consequence of the above equality we obtain  
\beq
\begin{split}
\upssone\left(\sci_{\sone}(\{q_0\},-\nabla\Psi_{\lambda_{\pm}})\right)&= \upssone\left(\sci_{\sone}(\{0\},-\nabla \widetilde{\Psi}_{\lambda_{\pm}})\right)=
\\ &= \upssone\left(\sci_{\sone}(\{0\}, -\nabla \varphi_{\lambda_{\pm}})\right) \star \upssone\left(\sci_{\sone}(\{0\},-{\nabla^2 \widetilde{\Psi}_{\lambda_{\pm}}(0)}_{\mid \cR})\right)=
\\  & = (-1)^{\dim \cN} \cdot \upssone\left(\sci_{\sone}\left(\{0\},-{\nabla^2 \widetilde{\Psi}_{\lambda_{\pm}}}_{\mid \cR} \right) \right)= 
\\ & = (-1)^{\dim \cN} \cdot \upssone\left(\sci_{\sone}\left(\{q_0\},-\nabla \Pi_{\lambda_{\pm}} \right) \right),
\end{split}
\eeq
which completes the proof.
\end{proof}

Our goal is to prove that  $\upssone\left(\sci_{\sone}\left(\{q_0\},-\nabla\Pi_{\lambda_+}\right)\right)\neq \upssone\left(\sci_{\sone}\left(\{q_0\},-\nabla\Pi_{\lambda_-}\right)\right)$. 
Therefore, following the construction of the equivariant Conley index given in Subsection \ref{conleyIndexSection}, we define the finite dimensional spaces 
$$\ds \bH^n= \left(T_{q_0}^{\perp} \Gamma(q_0) \oplus \bigoplus_{k=1}^{n} \bH_k\right)\cap \cR = \left(T_{q_0}^{\perp} \Gamma(q_0)\ominus \cN \right) \oplus \bigoplus_{k=1}^{n} \bH_k\subset \cR \subset \bH, \quad n\geq 1.$$
For simplicity of notation put   $\Pi_{\lpm}^n:=\Pi_{\lpm\mid\bH^n}$. We follow a similar reasoning as given in \cite{PCHRYST}.

\bl\label{main-lemma}
There exists $n_0 \in \bN$ such that  for any $n \geq n_0$ 
\begin{equation*}
\begin{split}
\chi_{\sone}(\ci_{\sone}(\{q_0\},-\nabla\Pi^n_{\lambda_-})) &= \chi_{\sone}(\ci_{\sone}(\{q_0\},-\nabla\Pi^{n_0}_{\lambda_-})) \\ 
\neq\chi_{\sone}(\ci_{\sone}(\{q_0\},-\nabla\Pi^{n_0}_{\lambda_+}))&= \chi_{\sone}(\ci_{\sone}(\{q_0\},-\nabla\Pi^n_{\lambda_+})).
\end{split}
\end{equation*} 
Moreover,  $\upssone\left(\sci_{\sone}\left(\{q_0\},-\nabla\Pi_{\lambda_+}\right)\right)\neq \upssone\left(\sci_{\sone}\left(\{q_0\},-\nabla\Pi_{\lambda_-}\right)\right)$.
\el
\begin{proof}
By assumption of Theorem \ref{main-theo}, for every $j=1,\ldots,j_0-1$ one can choose $k_j\in\bN$ such that $k_j^2<(\beta_j/\beta_{j_0})^2<(k_j+1)^2.$ Note that $k_1 \geq k_2 \geq \ldots \geq k_{j_0-2} \geq k_{j_0-1}.$ 
Taking into account that $\frac{(k_1+1)^2}{\beta_1^2} > \frac{1}{\beta_{j_0}^2},$  $\lambda_+=\frac{1+\varepsilon}{\beta_{j_0}}$ and that $\varepsilon$  is arbitrarily small, for fixed $n_0\geq k_1+1$ and $j=1,\ldots,m$ we obtain
\beq \label{szac}
n_0^2-\lambda_{\pm}^2\beta_j^2\geq n_0^2-\lambda_{\pm}^2\beta_1^2\geq n_0^2-\lambda_+^2\beta_1^2\geq \beta_1^2 \left(\frac{(k_1+1)^2}{\beta_1^2}-\lambda_{+}^2\right)>0.
\eeq
From the above formula and  equation \eqref{fourier} we get that for any $n \geq n_0$ the following equality holds 
$$
\morse(\nabla^2\Pi_{\lpm}^n)=\morse(\nabla^2\Pi_{\lpm}^{n_0}),
$$ 
where $\morse(\cdot)$ is the Morse index. Since $\Pi_{\lpm}^{n_0}$ is an isomorphism, for any $n \geq n_0$ we obtain 
$\ci_{\sone}(\{q_0\},-\nabla\Pi_{\lpm}^n)= \ci_{\sone}(\{q_0\},-\nabla\Pi_{\lpm}^{n_0}).$

Therefore  the  $\sone$-equivariant Conley index of the isolated   invariant set $\{q_0\}$ under the linear vector field $-\nabla\Pi_{\lpm}$ is the $\sone$-homotopy type of the spectrum $(E_{n,\pm})_{n=n_0}^{\infty},$ see \cite{IZYDOREK}, where $E_{n,\pm}$ is the same   pointed topological $\sone$-space for every $n\geq n_0.$ Consequently, the study of a change of the $\sone$-equivariant Conley index of $\{q_0\}$ under the linear vector field $-\nabla\Pi_{\lpm}$ is equivalent to the study of finite-dimensional $\sone$-equivariant Conley index $E_{n_0,\pm}=\ci_{\sone}(\{q_0\},-\nabla\Pi_{\lpm}^{n_0}).$

Note that since $[\lambda_-,\lambda_+] \cap \Lambda=\{1\slash \beta_{j_0}\}$, for $\lambda\in [\lambda_-,\lambda_+],\,k=1,\ldots,n_0$ and $j=1,\ldots,m$ we have
\beq\label{prop1}
k^2-\lambda^2 \beta_j^2 = 0\quad \text{ iff } \quad k=1,\,  j=j_0.
\eeq
Moreover,
\begin{equation}\label{prop4}
(1-\lambda^2_- \beta_{j_0}^2)(1-\lambda^2_+ \beta_{j_0}^2)=-\varepsilon^2(4-\varepsilon^2)<0.
\end{equation}
By formulas \eqref{prop1}, \eqref{prop4} and equation \eqref{fourier} we obtain that the only change of sign of eigenvalues of $\nabla^2\Pi_{\lambda}^{n_0}$ between $\lambda_-$ and $\lambda_+$ holds on $\bH_1$ and
\begin{equation}\label{prop5}
\morse\left(\nabla^2 {\Pi_{\lambda_-}}_{\mid \bH_1}\right) \neq \morse\left(\nabla^2 {\Pi_{\lambda_+}}_{\mid \bH_1}\right).
\end{equation}
Therefore the spectral decomposition of $\bH^{n_0}$ given by the isomorphism $-\nabla^2 \Pi_{\lpm}$ has the form
\beq
\bH^{n_0}=\bH_1\oplus\bW=\left(\bH_{1,\pm}^- \oplus \bH_{1,\pm}^+\right) \oplus \left(\bW^- \oplus \bW^+\right),
\eeq
where $\bW:=\left((T_{q_0}^{\perp} \Gamma(q_0) \ominus \cN)\oplus \bigoplus_{k=2}^{n_0} \bH_{k} \right),\,\bW^+,\,\bW^-$ don't depend on $\lambda_{\pm}$ and  
\beq\label{dimensions}
r_-:=\dim \bH^+_{1,-} \neq \dim \bH^+_{1,+} =: r_+.
\eeq
Hence $\ds \ci_{\sone} (\{q_0\},-\nabla \Pi_{\lpm}^{n_0})= S^{\bH_{1,\pm}^+} \wedge S^{\bW^+}$ and
\begin{equation}\label{prop6} 
\chi_{\sone} \left( \ci_{\sone} (\{q_0\},-\nabla \Pi_{\lpm}^{n_0})\right)= \chi_{\sone}\left(S^{\bH_{1,{\pm}}^{\pm}}\right) \star  \chi_{\sone}\left( S^{\bW^+}\right) \in U(\sone).
\end{equation}
By Remark \ref{propOfChi}  the element $\chi_{\sone}(S^{\bW^+})$ is invertible in the Euler ring  $U(\sone)$, hence it is sufficient to show
\begin{equation}\label{final formula}
\chi_{\sone}\left(S^{\bH_{1,-}^+}\right)\neq\chi_{\sone}\left(S^{\bH_{1,+}^+}\right).
\end{equation}
However, since $\bH_1= \lin \{e_i \cos t, e_i \sin t: i=1,\ldots,n\}$  and the action of  $S^1$ on $\bH_1$ is given by shift in time, the spaces $\bH^+_{1,\pm}$  are representations of   $\sone$ such that $\bH^+_{\pm} \approx_{\sone} \bR[r_{\pm},1].$ Hence by formula \eqref{ches} we obtain 
\[
\chi_{\sone}\left(S^{\bH^+_{1,\pm}}\right)= \chi_{\sone}\left(S^{\bR[r_{\pm},1]}\right)=    \bI - r_{\pm} \chi_{\sone}\left({\sone \slash \bZ_1}^+\right) \in U(\sone).
\]
Taking into account inequality \eqref{dimensions} we obtain \eqref{final formula} and consequently
\[
\chi_{\sone}\left(\ci_{\sone} (\{q_0\},-\nabla {\Pi_{\lambda_-}^{n_0}})\right) \neq \chi_{\sone}\left(\ci_{\sone} (\{q_0\},-\nabla {\Pi_{\lambda_+}^{n_0}})\right).
\]

Finally,  based on the definitions of $\upssone$ and  the equivariant Conley index and Remark \ref{stabilization of upsilon}, we obtain
\[
\upssone\left(\sci_{\sone}\left(\{q_0\},-\nabla\Pi_{\lambda_-}\right)\right)\neq \upssone\left(\sci_{\sone}\left(\{q_0\},-\nabla\Pi_{\lambda_+}\right)\right),
\]
because for $\lambda_+$ and $\lambda_-$ stabilization of the equivariant Conley index begins on the same level $n_0$. 
\end{proof}

\subsection{Proof of Theorem \ref{main-theo}}
\numsubsec 
\label{proof}
We are now in a position to prove the Symmetric Liapunov center theorem for minimal orbit.
\begin{proof}
To complete the proof of Theorem \ref{main-theo}, we will show  a change of the equivariant Conley index along the family $\cT = \Gamma(q_0) \times (0,+\infty)$ of trivial solutions which implies the existence of local bifurcation of non-stationary periodic solutions of system \eqref{newsys}, see Theorem \ref{sufficient}. Therefore, our aim is to prove that
\beq \label{finish}
\sci_G \left(G(q_0),-\nabla\Phi(\cdot,\lambda_-)\right)\neq \sci_G \left(G(q_0),-\nabla\Phi(\cdot,\lambda_+)\right),
\eeq
where $\lpm=\frac{1\pm\varepsilon}{\beta_{j_0}}$ and $\varepsilon>0$ is sufficiently small.

Combining  Lemmas \ref{main-lemma0} and \ref{main-lemma} we obtain 
\beq\label{eq1-proof}
\upssone\left(\sci_{\sone}(\{q_0\},-\nabla\Psi_{\lambda_+}\right)\neq \upssone\left(\sci_{\sone}(\{q_0\},-\nabla\Psi_{\lambda_-})\right),
\eeq
where $\Psi_{\lpm}=\Phi(\cdot,\lpm)_{\mid \bH} : \bH \to \bR$ and $\h=T_{q_0}^{\perp} G(q_0) \subset \h^1_{2\pi}$.

Since the pair $(\Gamma\times\sone,\{e\}\times\sone)$ is admissible and both $\nabla\Phi(\cdot,\lambda_-)$, $\nabla\Phi(\cdot,\lambda_+)$ are in the form of a compact perturbation of the identity, we can apply Theorem \ref{smash theorem for sci} together with equation \eqref{eq1-proof} to get
\beq\label{eq2-proof}
\upsg\left(\sci_G\left(G(q_0),-\nabla\Phi(\cdot,\lambda_-)\right)\right)\neq \upsg\left(\sci_G\left(G(q_0),-\nabla\Phi(\cdot,\lambda_+)\right)\right).
\eeq
Inequality  \eqref{eq2-proof} in a natural way implies inequality \eqref{finish}.
 
Since we have  proved a change of the equivariant Conley index on the segment $[\lambda_-,\lambda_+]$ and in this segment the only point where bifurcation can occur is $1/\beta_{j_0}$ we have just proved a local  bifurcation of non-stationary solutions of  problem \eqref{sys} from the orbit $G(q_0)\times\{1/\beta_{j_0}\}$. As we know, they correspond to $2\pi\lambda (\approx \frac{2\pi}{\beta_{j_0}})$-periodic solutions of the system \eqref{newsys} in a neighborhood of the orbit $G(q_0)$. From the assumptions  it follows that $1/(k\cdot \beta_{j_0})\notin \Lambda$ holds true  for $k\geq 2$  (i.e. the orbit $G(q_0)\times\{\frac{1}{k\cdot \beta_{j_0}}\}$ can not be an orbit of bifurcation), therefore the bifurcating orbits have the minimal period close to $\frac{2\pi}{\beta_{j_0}}$, which completes the proof.

\end{proof}

\br
It was shown in \cite{GORY1} that $\upsg\left(\sci_{G}\left(G(q_0),-\nabla\Phi\right)\right)$ is equal to degree for $G$-equivariant gradient maps $\nabla_G\textrm{-}\mathrm{deg}(\nabla \Phi, \cO)$ defined in \cite{GORY}, where $\cO\subset\bH^1_{2\pi}$ is a $G$-invariant subset such that $\nabla \Phi^{-1}(0) \cap \cO = G(q_0)$. Moreover, the change of this degree along the family of trivial solutions implies the bifurcation of a connected set of solutions of equation $\nabla\Phi(q,\lambda)=0$.

 Since we have proved that $\upsg\left(\sci_{G}\left(G(q_0),-\nabla\Phi(\cdot,\lambda_-)\right)\right)\neq \upsg\left(\sci_{G}\left(G(q_0),-\nabla\Phi(\cdot,\lambda_+)\right)\right)$,
there exists a connected family $\cC$ of non-stationary periodic solutions of equation \eqref{newsys} emanating from the orbit $\Gamma(q_0)$ which satisfies the Rabinowitz-type alternative.
\er

%%%%%%%%%%%%%%%%%%%%%%%%%%%%%%%%%%%%%%%%%%%%%%%%%%%%%%%

\section{Applications} \label{applications}
\numsec

In order to show the strength of our main result i.e. the Symmetric Liapunov center theorem for minimal orbit, see  Theorem \ref{main-theo}, in this section we apply it to the study of periodic solutions of the Lennard-Jones and Schwarzschild $2$-, $3$-body problems with $\Gamma=\sotwo$.

 Consider  $(\bR^2)^N$ as an orthogonal representation of $\sotwo$ with $\sotwo$-action defined by 
$\sotwo \times (\bR^2)^N \ni (\gamma,(q_1,\ldots,q_N)) \to (\gamma q_1,\ldots, \gamma q_N) \in (\bR^2)^N.$
Define an open $\sotwo$-invariant subset $\Omega=\{q = (q_1,\ldots,q_N) \in (\bR^2)^N : q_i \neq q_i \text{ for } i \neq j\}$ and note that if $q_0 \in \Omega$ then the isotropy group $\sotwo_{q_0}$ is trivial.
Recall that $\sigma(S)$ stands for the spectrum of a symmetric matrix $S$ and denote by $\m(\alpha)$ the multiplicity of an eigenvalue $\alpha \in \sigma(S).$

\subsection{The Lennard-Jones problem}
\numsubsec

The Lennard-Jones potential is used to model the nature and stability of small clusters of interacting particles in crystal growth, random geometry of liquids, and in the theory of homogeneous nucleation, see for instance \cite{WD}. The potential also appears in molecular dynamics to simulate many particles systems ranging from solids, liquids, gases and biomolecules of Earth.

It is well known in chemistry and chemical physics that the stability of some molecular structures are closely related with the local minima of the corresponding potential. Also in the analysis of the native structure of a protein, it is necessary to find the lowest energy configuration of a molecular system. In general when it is possible to find the global minimum of a potential energy surface, we can get a global optimization of the problem, saving  money and laboratory time in this way. Unfortunately, this is a difficult task in general, and the researchers on the subject develop global optimization methods on simpler systems, one of the most useful in this direction is the Lennard-Jones potential, which has been used in the analysis of clusters in nanomaterials  in the last times, see for instance \cite{WD} and the references therein.

The Lennard-Jones potential is given by 
$$U = \varepsilon \sum \left[ \left( \frac{\sigma}{r_{ij}} \right)^{12} - 2\left( \frac{\sigma}{r_{ij}} \right)^6 \right],$$
where $\varepsilon$ represents the minimum value of the potential energy and $\sigma$ is the minimum distance from the origin on the $x$-axis (when the potential is repulsive). That is, all molecules are attracting each other when they are close enough, the intensity of this attraction force decreases when the molecular distance increases.  

For simplicity we assume that  $\varepsilon = \sigma = 1.$ So we
consider $N$-particles with equal mass $m$ moving in the $2$-dimensional Euclidean space. The forces between two particles are given by the Lennard-Jones potential. Let $q_i$ denotes the position of the  $i$-th particle in an inertial coordinate system and let $q=(q_1,\ldots,q_N) \in (\bR^2)^N.$ Choosing the units of mass, length and time conveniently one can define  the Lennard-Jones potential $U : \Omega \to \bR$ as
\beq
\label{ljpot}
U(q)=\sum_{1 \leq  i < j \leq N} \left(\frac{1}{\mid q_i-q_j\mid^{12}}  - \frac{2}{\mid q_i-q_j\mid^{6}} \right).
\eeq

Note that $U : \Omega \to \bR$ is smooth and $\sotwo$-invariant. The Lennard-Jones problem has been widely studied in \cite{CLC1}, where the authors show the existence of families of periodic orbits, where the mutual distances among the particles remain constant along the motion, i.e. the particles behave as a rigid body. These especial kind of periodic orbits are called {\it relative equilibria}.
The aim of this section is to show the strength of our main Theorem \ref{main-theo}, by showing 
new families of  periodic solutions in the Lennard-Jones problem  defined by equation \eqref{ljpot}.

\subsubsection{Case $N=2$}
\numsubsubsec
First consider   Lennard-Jones $2$-body problem. In this case is easy to verify that the following equality holds $(\nabla U)^{-1}(0) \cap \Omega=\{(q_1,q_2) \in \Omega : q_1=-q_2 \text{ and } \mid q_1-q_2\mid=1\},$ see Theorem 1 of \cite{CLC1}. In other words $(\nabla U)^{-1}(0) \cap \Omega = \Gamma(q_0),$ where $q_0=(0,1/2,0,-1/2)$ i.e. the set of critical points of the potential $U$ consists of one orbit $\Gamma(q_0).$  Moreover, it was proved in \cite{CLC1}
that the orbit $\Gamma(q_0)$ consists of minima of $U$ and $U(\Gamma(q_0))=-1.$ Since the action of $\sotwo$ on $\Omega$ is free, the isotropy group $\sotwo_{q_0}$ is trivial. It is easy to check that $\sigma(\nabla^2 U(q_0))=\{0,144\}$ and $mult(0)=3$, $mult(144) = 1.$ 

 We have just shown that all assumptions of  Theorem  \ref{main-theo} are fulfilled with $j_0=1$ and $\beta_1=12$. Applying this theorem we obtain the existence of non-stationary periodic solutions of system \eqref{newsys} in any neighborhood of the orbit $\sotwo(q_0).$ Moreover, the minimal period of these solutions are close to $\pi \slash 6 \approx 0,5235.$ 

\br  The existence of a family $\cF$ of non-stationary periodic solutions of system \eqref{newsys}  (relative equilibria), has been proved in Theorem 3 of \cite{CLC1}. In any neighborhood of the orbit $\sotwo(q_0)$ there is a member of this family. Moreover, on that paper the authors proved that the  lower estimation of minimal period of members of $\cF$ is equals to $\frac{7 \pi}{6}\left( \frac{7}{32}\right)^{\frac{1}{6}} \approx 2.8450$, see Theorem 3 of \cite{CLC1}.  Since $\frac{7 \pi}{6}\left( \frac{7}{32}\right)^{\frac{1}{6}} > \pi \slash 6$, we  point out that the non-stationary periodic solutions of \eqref{newsys} with the Lennard-Jones potential,  whose existence we have proved in a sufficiently small neighborhood of $\sotwo(q_0)$ are different from the relative equilibria whose existence was proved in \cite{CLC1}.\er

\subsubsection{Case $N=3$}
\numsubsubsec

For the Lennard-Jones $3$-body problem.  Set
$$
\ba{ccl} 
q_{01} & = & (a/2,0,0,0,-a/2,0),\\
q_{02} & = & (0,0,a/2,0,-a/2,0),\\
q_{03} & = & (a/2 ,0 ,-a/2,0 , 0,0), \\
q_{04} & = & (1 / \sqrt{3}) (1,0,\cos \alpha,\sin \alpha,\cos \beta,\sin \beta),\\
q_{05} & = & (1 / \sqrt{3})  (1,0,\cos \beta,\sin \beta,\cos \alpha,\sin \alpha),
\ea
$$
where $a=\left(\frac{2731}{43}\right)^{1/6}, \alpha = 2 \pi / 3$ and $ \beta=4 \pi/3.$

Put $\Gamma=\sotwo$. We know that the following equality holds $\ds (\nabla U)^{-1}(0) = \bigcup_{i=1}^5 \Gamma(q_{0i}),$  see Theorem 6 of \cite{CLC1}. On that paper the authors also  proved that the orbits of critical points $\Gamma(q_{04}), \Gamma(q_{05})$ are minima of the energy potential $U$ given by \eqref{ljpot} and $U(\Gamma(q_{04}))=U(\Gamma(q_{05}))=-3$. Moreover,  the isotropy groups $\Gamma_{q_{04}},\Gamma_{q_{05}}$ are trivial. It is not difficult to verify that $\sigma(\nabla^2 U(q_{04})) = \sigma(\nabla^2 U(q_{05}))= \{0,108,216\}, $ and that $\m(0)=3, \m(108)=2, \m(216)=1.$ Note that all assumptions of Theorem \ref{main-theo} are fulfilled at orbits $\Gamma(q_{04}), \Gamma(q_{05})$ with $\beta_1 = 6 \sqrt{6}, \beta_2 = 6 \sqrt{3} $ and $j_0=1$ or $j_0=2$. Applying this theorem we obtain the existence 
of periodic solutions in any  neighborhood of the orbits  $\Gamma(q_{04}),$ $ \Gamma(q_{05})$  with minimal periods close to $\pi /(3 \sqrt{3})$ and $\pi /(3 \sqrt{6}).$

\br It is known that $0 \in \sigma(\nabla^2 U(q_{01})) =\sigma(\nabla^2 U(q_{02})) = \sigma(\nabla^2 U(q_{03})), \m(0)=3$ and that $\sigma(\nabla^2 U(q_{01})) \cap (0,+\infty) \neq \emptyset,$ see Section 4.1.1 of \cite{CLC1}. That is why the orbits $\Gamma(q_{01}),$ $ \Gamma(q_{02}),$ $ \Gamma(q_{03})$ are isolated and degenerate. We underline that they are not orbits of minima of the potential $U.$ That is why one can not apply Theorem \ref{main-theo} to the study of existence of non-stationary periodic solutions of system \eqref{newsys} in a neighborhood of these orbits. At this moment this problem is far from being solved.
\er

\subsection{The Schwarzschild 3-body problem}
\numsubsec 

The potential in which we are interested comes from relativistic physics. It was introduced in 1916 by  Schwarzschild \cite{SCH} in order to give a solution to Einstein's equations for the gravitational field of an uncharged spherical non-rotating mass, which trough a classical formalism provides the Binet-type equations. Again, as in the Lennard-Jones problem the corresponding force (after a normalization of coordinates using cosmological background), is given by the minus  gradient of the so called Schwarzschild potential, which
 has the simple form
\begin{equation}\label{schwarz-pot}
U = \frac{A}{r} + \frac{B}{r^{3}}.
\end{equation}

 The Schwarzschild potential was tackled into the framework of dynamical systems and celestial mechanics by  Mioc and collaborators, see for instance \cite{MPS}, \cite{STMI} and the references therein. This new and original approach to study the dynamics of particles moving under this potential has been very useful in astrophysics for the analysis of theoretical black holes or the motion of a galaxy far enough that you can consider it as a single object. It has also been used in cosmology for the analysis of clusters of galaxies.  The case $A<0<B$ which concerns with the main results of this paper models the photogravitational field of the Sun, see \cite{APCS} and the references therein for more details.

For  $1 \leqq i < j \leqq 3$  choose $A_{ij} < 0 <B_{ij},$ define $\cU_{ij} \in  C^2((0,+,\infty),\bR)$ by $\ds \cU_{ij}(r_{ij})= \frac{A_{ij}}{r_{ij}} + \frac{B_{ij}}{r_{ij}^3}$ and  $\cU \in C^2((0,+\infty)^3,\bR)$ by $\ds \cU(r_{12},r_{13},r_{23})=\sum_{1 \leqq i < j \leqq 3} \cU_{ij}(r_{ij}).$ 

 The Schwarzschild $3$-body potential $U : \Omega \to \bR$ is defined by
\beq 
\label{schpot}
U(p)=\cU(r_{12}(q),r_{13}(q),r_{23}(q)) = \sum_{1 \leqq i < j \leqq 3} \cU_{ij}(r_{ij}(q)),
\eeq 
where $r_{ij}(q)=\|q_i-q_j\|.$ We observe that the Schwarzschild potential is smooth and $\sotwo$-invariant. As for the Lennard-Jones problem, our goal here is to show the strength of our main Theorem \ref{main-theo} to show the existence of new families of  periodic solutions of system \eqref{newsys}  with the Schwarzschild potential $U$ defined by formula \eqref{schpot}. 

In the following lemma we describe local non-degenerate minima of the potentials  $\cU_{ij}$, for $1 \leqq i < j \leqq 3.$

\bl  \label{sch1}  For $1 \leqq i < j \leqq 3$   define $\ds r_{ij}^0=\sqrt{-\frac{3B_{ij}}{A_{ij}}}$ and note that 
\be 
\item $\ds \cU'_{ij}(r)=0$ if and only if $\ds r=r_{ij}^0,$
\item $\ds \cU''_{ij}(r_{ij}^0) =\frac{3B_{ij}}{(r_{ij}^0)^5} > 0.$
\ee
\el

Now we characterize the critical points of the potential $U.$

\bl \label{sch2} Fix $q_0=(q_{10},q_{20}, q_{30}) \in \Omega$ such that vectors $q_{10}, q_{20}, q_{30}$ are not collinear. Then $q_0$ is  a critical point of $U$ and  iff
$\ds r_{12}(q_0)=\sqrt{-\frac{3B_{12}}{A_{12}}}, r_{13}(q_0)=\sqrt{-\frac{3B_{13}}{A_{13}}}, r_{23}(q_0)=\sqrt{-\frac{3B_{23}}{A_{23}}}.$
\el 
\begin{proof}
First, we observe  that  
\beq \label{sysu}
\ba{rcl}
\frac{\partial U}{\partial q_1}(q) & = &\ds \frac{1}{r_{12}(q)} \cU'_{12}(r_{12}(q))   (q_1-q_2) + \frac{1}{r_{13}(q)} \cU'_{13}(r_{13}(q))   (q_1-q_3),\\
\frac{\partial U}{\partial q_2}(q) & = &\ds   \frac{1}{r_{12}(q)} \cU'_{12}(r_{12}(q))   (q_2-q_1) +  \frac{1}{r_{23}(q)} \cU'_{23}(r_{23}(q))    (q_2-q_3), \\
\frac{\partial U}{\partial q_3}(q) & = & \ds \frac{1}{r_{13}(q)} \cU'_{13}(r_{13}(q))   (q_3-q_1) + \frac{1}{r_{23}(q)} \cU'_{23}(r_{23}(q))   (q_3-q_2).
\ea
\eeq
Since   every equation  in the right hand side of system \eqref{sysu} is a linear combination of  linearly independent vectors,   applying Lemma  \ref{sch1} we obtain
$$\nabla U(q_0)=0 \quad \text{ iff } \quad \cU'_{12}(r_{12}(q_0))=0, \quad \cU'_{13}(r_{13}(q_0))=0, \quad \cU'_{23}(r_{23}(q_0))=0 $$ 
$$
\text{ iff } \quad r_{12}(q_0)=\sqrt{-\frac{3B_{12}}{A_{12}}}, r_{13}(q_0)=\sqrt{-\frac{3B_{13}}{A_{13}}}, r_{23}(q_0)=\sqrt{-\frac{3B_{23}}{A_{23}}},$$
which completes the proof.
\end{proof}

 Fix $q_0=(q_{10},q_{20}, q_{30}) \in \Omega$ such that vectors $q_{10}, q_{20}, q_{30}$ are not collinear and choose constants $A_{12}, A_{13}, A_{23} < 0 < B_{12}, B_{13}, B_{23}$ such that $$ \|q_{10}-q_{20}\| =\sqrt{-\frac{3B_{12}}{A_{12}}}, \quad \|q_{10}-q_{30}\|=\sqrt{-\frac{3B_{13}}{A_{13}}}, \quad \|q_{20}-q_{30} \|=\sqrt{-\frac{3B_{23}}{A_{23}}}.$$

From Lemma \ref{sch2} it follows that $\nabla U(q_0)=0.$ Since the Schwarzschild potential $U$ is $\sotwo$-invariant, $\sotwo(q_0) \subset (\nabla U)^{-1}(0).$ Additionally, from Lemma \ref{sch1} we obtain that $\sotwo(q_0)$  is an isolated minimal orbit of critical points of $U.$ Finally, since the orbit $\sotwo(q_0)$ is minimal, all the nonzero eigenvalues of $\nabla^2 U(q_0)$ are positive i.e. $\sigma(\nabla^2 U(q_0))=\{0\} \cup \left(\sigma(\nabla^2 U(q_0)) \cap (0,+ \infty)\right)$ $ = \{0, \beta_1^2,\ldots,\beta_m^2\}.$ Fix $\beta_{j_0}$ such that $\beta_j \slash \beta_{j_0}\not \in \bN$ for $j \neq j_0.$ 

We have shown that all assumptions of Theorem \ref{main-theo} are fulfilled.  Applying this theorem we obtain the existence 
of an orbit of periodic solutions of system \eqref{newsys} in any  neighborhood of the orbit  $\sotwo(q_0)$  with minimal period close to $2\pi / \beta_{j_0}.$

 In order to illustrate the above reasoning, we apply it to a simple example.
 
\bex 
Set $A_{12}=- 3 \slash 2, A_{13}=-1, A_{23}=-3 \slash 5, B_{12}=1 \slash 2, B_{13} = 1 \slash 3, B_{23}=1 \slash 5,$  $q_{10}=(\sqrt{2} \slash 2,0), q_{20}=(0,\sqrt{2} \slash 2), q_{30}=((\sqrt{2}+\sqrt{6}) \slash 4,(\sqrt{2}+\sqrt{6}) \slash 4)$ and $q_0=(q_{10},q_{20},q_{30}).$ It is easy to verify that 
$$ \|q_{10}-q_{20}\| =1=\sqrt{-\frac{3B_{12}}{A_{12}}}, \quad \|q_{10}-q_{30}\|=1=\sqrt{-\frac{3B_{13}}{A_{13}}}, \quad \|q_{20}-q_{30} \|=1=\sqrt{-\frac{3B_{23}}{A_{23}}}.$$

It shows that $\sotwo(q_0) \subset (\nabla U)^{-1}(0).$ The characteristic polynomial of $\nabla^2 U(q_0)$ is equal to $w(x)=\frac{1}{5} x^3(6 x^3-62 x^2+225 x-243)$ and its roots are   $\beta_0=0,\beta_1^2 \approx 2.027 ,$ $ \beta_2^2 \approx 3.475 , \beta_3^2 \approx 6.897$. Moreover, $\m(\beta_0)=3, \m(\beta_1)=\m(\beta_2)=\m(\beta_3)=1.$ Since $\beta_i \slash \beta_j \notin \bN$
for $i \neq j,$   in any  neighborhood of orbit  $\sotwo(q_0)$  there are orbits of periodic solutions of system  \eqref{newsys} with minimal periods close to $2\pi / \beta_1, 2 \pi \slash \beta_2, 2 \pi \slash \beta_3.$
\eex

\subsection*{Acknowledgements} 
The first author (EPC) has been partially supported by {\it Asociaci\'on Mexicana de Cultura A.C.}.

\bibliographystyle{abbrv}
\bibliography{20170910bibliography}

 \end{document}